\documentclass[11pt,twoside]{article} 
\usepackage{epsfig,color}
\usepackage{amsmath}
\usepackage{amssymb, geometry,url}

\newcommand{\diag}{\mathop{\textrm{diag}}}
\newcommand{\bA}{\boldsymbol{A}}
\newcommand{\bSigma}{\boldsymbol{\Sigma}}
\newcommand{\bZ}{\boldsymbol{Z}}
\newcommand{\bW}{\boldsymbol{W}}
\newcommand{\bx}{\boldsymbol{x}}
\newcommand{\bd}{\boldsymbol{d}}
\newcommand{\bz}{\boldsymbol{z}}
\newcommand{\blambda}{\boldsymbol{\lambda}}
\newcommand{\bX}{\boldsymbol{X}}
\newcommand{\bone}{\boldsymbol{1}}

\newcommand{\ubar}[1]{\underaccent{\bar}{#1}}
\usepackage{accents}
\usepackage{listings}
\usepackage[dvipsnames]{xcolor}

\lstdefinelanguage{julia}
{
keywordsprefix=\@,
morekeywords={
exit,whos,edit,load,is,isa,isequal,typeof,tuple,ntuple,uid,hash,finalizer,convert,promote,
subtype,typemin,typemax,realmin,realmax,sizeof,eps,promote_type,method_exists,applicable,
invoke,dlopen,dlsym,system,error,throw,assert,new,Inf,Nan,pi,im,begin,while,for,in,return,
break,continue,macro,quote,let,if,elseif,else,try,catch,end,bitstype,ccall,do,using,module,
import,export,importall,baremodule,immutable,local,global,const,Bool,Int,Int8,Int16,Int32,
Int64,Uint,Uint8,Uint16,Uint32,Uint64,Float32,Float64,Complex64,Complex128,Any,Nothing,None,
function,type,typealias,abstract,get_node,add_edge,create_estimation_model,set_solution,
solve, get_solution, solve_ss_problem, create_estimation_problem, addnode
},
morekeywords = [2]{triggered_by,compute_time,trigger_during_busy,send_on,delay,send_wait,start},
sensitive=true,
morecomment=[l]{\#},
morestring=[b]',
morestring=[b]"
}
\usepackage{multirow}

\usepackage[colorlinks=true,linkcolor=blue,citecolor=blue,urlcolor=blue]{hyperref}

\usepackage{fancyhdr}
\begin{document}

\title{A Julia Framework for Graph-Structured\\ Nonlinear Optimization}

\author{David L. Cole${}^{\dag}$, Sungho Shin${}^{\ddag}$, and Victor Zavala${}^{\dag\ddag}$\thanks{Corresponding Author: victor.zavala@wisc.edu}
}

\date{\small
  ${}^\dag$Department of Chemical and Biological Engineering, \\[0in]
  University of Wisconsin-Madison, Madison, WI\\[.05in]
  ${}^\ddag$Mathematics and Computer Science Division, \\
  Argonne National Laboratory, Argonne, IL\\[-2in]
}

\maketitle

\begin{abstract}
Graph theory provides a convenient framework for modeling and solving structured optimization problems. Under this framework, the modeler can arrange/assemble the components of an optimization model (variables, constraints, objective functions, and data) within nodes and edges of a graph, and this representation can be used to visualize, manipulate, and solve the problem. In this work, we present a {\tt Julia} framework for modeling and solving graph-structured nonlinear optimization problems.  Our framework integrates the modeling package {\tt Plasmo.jl} (which facilitates the construction and manipulation of graph models) and the nonlinear optimization solver {\tt MadNLP.jl} (which provides capabilities for exploiting graph structures to accelerate solution). We illustrate with a simple example how model construction and manipulation can be performed in an intuitive manner using {\tt Plasmo.jl} and how the model structure can be exploited by {\tt MadNLP.jl}. We also demonstrate the scalability of the framework by targeting a large-scale, stochastic gas network problem that contains over 1.7 million variables. 
\end{abstract}


\section{Introduction}

Modeling and solving large nonlinear optimization problems is essential in diverse applications such as stochastic optimization, dynamic optimization, PDE-constrained optimization, and network optimization \cite{shin2021}. The complexity of such problems continuously pushes the boundary of existing computational tools and limits application scope. To overcome these challenges, it is necessary to develop tools that can facilitate the detection, manipulation, and exploitation of problem structure \cite{colombo2009structure,gondzio2009exploiting,gondzio2003parallel, grothey2009structure, jalving2021, wan2019, yoshio2021,zavala2008}. 
\\

Recently, it has been proposed to represent optimization problem structures in the form of graphs \cite{allman2019,berger2021gboml_tutorial, daoutidis2019decomposition, jalving2019, jalving2021, mitrai2020decomposition,mitrai2021blockmodeling,mitrai2021stochastic,tang2018optimal}. Under a graph representation, the components of an optimization problem (variables, constraints, and objectives) are assigned to nodes and edges  \cite{jalving2021,shin2020MadNLP}. Representing the problem structure as a graph has several benefits; specifically, the graph structure can be used to visualize and manipulate the model (e.g., graph partitioning) using powerful tools such as {\tt Metis} \cite{karypis1998} or {\tt KaHyPar} \cite{schlag2016}. Moreover, the graph structure can be communicated to optimization solvers and this facilitates the use of structure-exploiting linear algebra strategies inside nonlinear optimization solvers such as Schur decomposition \cite{bartlett2006qpschur, kang2014,laird2008large, laird2011parallel, rao1998, rodriguez2020,word2014efficient,zavala2008, zhu2009exploiting} and Schwarz decomposition  \cite{frommer2001,na2020,shin2020b,shin2020MadNLP,shin2020a}. 
\\

There are currently a few software frameworks that enable the construction and solution of structured optimization problems: {\tt Pyomo} (available in Python) \cite{hart2017pyomo, hart2011pyomo}, the Graph-Based Optimization Modeling Language, {\tt GBOML} \cite{berger2021gboml_tutorial,berger2021gboml}, and {\tt Plasmo.jl} (available in Julia) \cite{jalving2021}. {\tt Pyomo} uses a modeling extension called a “network” to model structured optimization problems; this approach creates ports (collections of variables) with arcs placed to link objects on separate ports. {\tt Pyomo} also has additional capabilities for working with the graph structure of stochastic optimization problems \cite{watson2012pysp}. {\tt GBOML} builds an optimization problem into blocks, where blocks are either nodes (containing variables, constraints, and/or objectives) or hyperedges (containing constraints) of a hypergraph. {\tt GBOML} is geared more specifically towards mixed-integer linear programs and is designed to facilitate efficient modeling and potentially enable decomposition schemes to exploit problem structure \cite{berger2021gboml_tutorial}. {\tt Plasmo.jl} places the problem components (variables, constraints, objectives, and data) into an abstract modeling object called an OptiGraph that is composed of OptiNodes (containing the modeling components) and OptiEdges (capturing the structural connectivity between components). These packages allow the modeler to define problem structure directly when constructing the model. This approach differs from approaches that aim to detect graph structures after the model has been built \cite{allman2019,daoutidis2019decomposition,mitrai2020decomposition,mitrai2021blockmodeling, mitrai2021stochastic,tang2018optimal}. The internal OptiNode models of {\tt Plasmo.jl} is constructed based on {\tt JuMP.jl} \cite{Dunning2017jump}, which is an algebraic modeling package available in Julia. {\tt JuMP.jl} provides the algebraic modeling user interface, which is inherited by {\tt Plasmo.jl}, and also provides the automatic differentiation capabilities. 
\\

Graph structures are exploited in nonlinear optimization solvers during the computation of Newton-like search steps (at the linear algebra level) using techniques such as Schur and Schwarz decomposition. Shin  and co-workers recently implemented a nonlinear, interior-point optimization solver in {\tt Julia} called {\tt MadNLP.jl} that provides these types of decomposition techniques \cite{shin2020MadNLP}. This solver also offers an extension called {\tt MadNLPGraph.jl}, which directly takes an OptiGraph object from {\tt Plasmo.jl} and exploits this structure during the evaluation of objective/constraint functions and of derivative information and during the Newton step computation. 
\\

Recent work by Jalving and co-workers showcased {\tt Plasmo.jl} capabilities, discussed the underlying OptiGraph abstraction, and provided an interface to the {\tt PIPS-NLP}  solver (implemented in C++) for exploiting graph structures using Schur decomposition \cite{chiang2014, jalving2021}. Recent work by Shin and co-workers \cite{shin2020MadNLP} introduced {\tt MadNLP.jl} and explored the use of Schwarz decomposition for exploiting graph structures provided by {\tt Plasmo.jl}. In this work, we provide an in-depth overview on capabilities that arise from the integration of {\tt Plasmo.jl} and {\tt MadNLP.jl} (see Figure \ref{fig:graphical_abstract}) and present Schur decomposition capabilities that have been implemented in {\tt MadNLP.jl}. Schur decomposition provides a flexible and robust decomposition paradigm for exploiting diverse structures communicated by {\tt Plasmo.jl} (e.g., time, space, or scenario structures). Moreover, we provide a detailed illustrative example to highlight diverse modeling and solution capabilities and we demonstrate scalability of the framework using a large-scale optimization problem that arises in the context of stochastic optimal control of natural gas networks and that contains over 1.7 million variables. All the code needed to reproduce the results of the paper can be found in \url{https://github.com/zavalab/JuliaBox/tree/master/GraphNLP}.
\\

The manuscript is structured as follows. Section 2 provides an overview of graph-based modeling and solution as well as implementations in {\tt Plasmo.jl} and {\tt MadNLP.jl}. Section 3 provides an illustrative example to highlight modeling, visualization, and partitioning capabilities of {\tt Plasmo.jl}. Section 4 provides a case study for a stochastic, nonlinear optimization problem arising in gas networks. Section 5 provides conclusions and a perspective on future work.

%%%%%%%%%%%%%%%%%%%%%%%%%%%%%%%%%%%%%%%%%%
\section{Overview of {\tt Plasmo.jl} and {\tt MadNLP.jl}}

\begin{figure}[!htp]
  \centering
  \includegraphics[scale=.5]{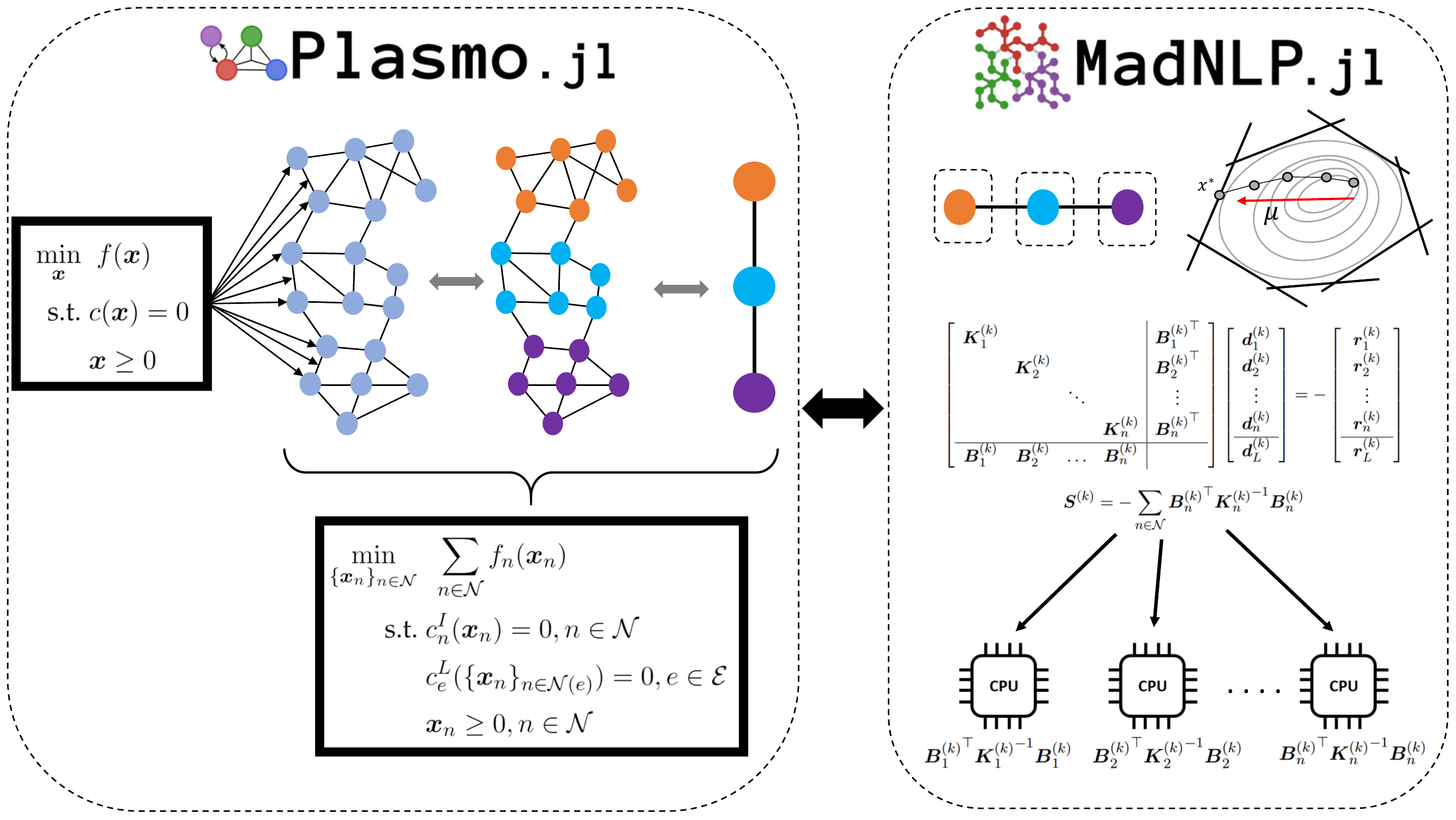}
  \caption{Visualization of the interaction of {\tt Plasmo.jl} and {\tt MadNLP.jl}. {\tt Plasmo.jl} enables a generic optimization problem to be modeled as a graph-based model, which can then be partitioned and/or aggregated. The graph-based model from {\tt Plasmo.jl} can be passed to the solver package {\tt MadNLP.jl}. {\tt MadNLP.jl} is an interior-point solver which exploits the problem structure that results from graph-based modeling to parallelize function and derivative evaluations and to solve the linear systems in parallel using Schur or Schwarz decompositions.}
  \label{fig:graphical_abstract}
\end{figure}

Graph-based modeling and solution of optimization problems have several benefits made possible by graph analysis tools. We use the packages {\tt Plasmo.jl} for constructing and partitioning problems as a graph-based model and use {\tt MadNLP.jl} for solving the graph-based models by exploiting the structure. In this section, we discuss the implementation of graph-based modeling tool {\tt Plasmo.jl} and structure-exploiting solver {\tt MadNLP.jl} and the interface between the packages. The interplay between {\tt Plasmo.jl} and {\tt MadNLP.jl} is illustrated in Figure \ref{fig:graphical_abstract}.

\subsection{{\tt Plasmo.jl}}

A compact representation of a graph-structured nonlinear optimization problem was proposed in \cite{jalving2021} and is given by:
\begin{subequations}\label{eq:graph_model}
\begin{align} 
\min_{\{\bx_n\}_{n \in {\mathcal{N}}}} &\;\;   \sum_{n \in \mathcal{N}} f_n(\bx_n)\\
\textrm{s.t.} &\; c_n^I(\bx_n) = 0,\quad n \in {\mathcal{N}}\\
&\; c_e^L(\{\bx_n\}_{n \in \mathcal{N}(e)})=0,\quad e \in \mathcal{E}\\
&\; \bx_n \geq 0,\quad n \in \mathcal{N}
\end{align} 
\end{subequations}
\noindent The graph $\mathcal{G}(\mathcal{N},\mathcal{E})$ has associated node set $\mathcal{N}$ and (undirected) edge set $\mathcal{E}$. The neighborhood of any node $n$ is denoted by $\mathcal{N}(n)$, while the set of nodes supporting an edge $e$ are represented by $\mathcal{N}(e)$. We also have that $\mathcal{E}(n)$ represents the set of all edges that are incident to node $n$. Symbol $\bx_n$ denotes the decision variables on the node $n$ and $f_n(\cdot)$ represents the objective function on node $n$. Internal equality constraint functions (i.e., constraints contained within a node) are denoted  by $c_n^I(\cdot)$. Linking constraint functions (i.e., constraints that embed variables in multiple nodes) are given by $c_e^L(\cdot)$. 
\\

The graph-based formulation \eqref{eq:graph_model} provides a unifying abstraction to capture structures appearing in diverse applications such as dynamic optimization (graph is a line), stochastic optimization (graph is a tree), and PDE optimization (graph is a discretization mesh).  In fact, most sparse nonlinear optimization problems of the form: 
\begin{subequations}\label{eq:optimization_model}
\begin{align} 
\min_{\bx} &\;\;   f(\bx)\\
\textrm{s.t.} &\; c(\bx) = 0\\
&\; \bx \geq 0
\end{align} 
\end{subequations}
can in principle be rearranged in the form \eqref{eq:graph_model} by using a proper allocation/assignment of the problem components into nodes and edges. This allocation can be done manually by the modeler (using domain-specific insight) or automatically using graph analysis tools (e.g., partitioning and community detection). Representing the problem as a graph can be particularly useful for problems that exhibit sparse structures and a high degree of modularity. 
\\

We should highlight that there is significant flexibility in how components can be assigned to nodes and edges in a graph abstraction; for instance, we could allocate a single variable/constraint to a node or we could embed the entire set of variables/constraints into a single node.  Such flexibility can be exploited for developing diverse solution strategies or to facilitate visualization. In principle, even dense/unstructured problems can be modeled as \eqref{eq:graph_model} form. However, for those problems, there is no practical benefit of using graph-based modeling packages and structure-exploiting solution algorithms.
\\

The {\tt Julia} package {\tt Plasmo.jl} has been developed to facilitate the representation of optimization problems  using graph structures. The implementation details of {\tt Plasmo.jl} are presented in  \cite{jalving2021}. {\tt Plasmo.jl} is based on a general modeling abstraction that is called an OptiGraph. An OptiGraph is an abstract object given by an undirected hypergraph that is composed of sets of OptiNodes and OptiEdges. An OptiNode is a self-contained optimization problem with components given by variables, constraints, objectives, and data; the OptiNode internally contains a {\tt JuMP.jl} model object. Connectivity across OptiNodes is captured in the form of OptiEdges, which are linking constraints that connect components of two (or more) OptiNodes (hence creating a hypergraph). The OptiGraph abstraction also enables hierarchical graph representations under which an OptiNode can be an OptiGraph itself; this facilitates the construction of problems with complex nested structures and manipulation of such structures (e.g., via aggregation and partitioning). In the following discussion, we will refer to OptiGraphs, OptiNodes, and OptiEdges as graphs, nodes, and edges. 
\\

Hierarchical graphs provide substantial modeling flexibility; for instance, the nodes within a graph can be partitions of variables, constraints, objectives, and/or data, and such partitions can be created on-the-fly using the {\tt partition} function. Hierarchical graphs also allow the user to place collections of nodes onto a subgraph, creating a secondary partition of the problem. While both nodes and subgraphs act as partitions, the nodes are defined in the initial problem formulation and do not change, while the subgraphs can be altered and partitioned further using functions within {\tt Plasmo.jl}. 
\\

Hierarchical graphs containing multiple levels of subgraphs can be aggregated/collapsed into a single node using the {\tt aggregate} function. The {\tt aggregate} function also allows modelers to convert the entire graph (including all nodes, edges, and subgraphs) into a single node containing the entire model components (a generic {\tt JuMP.jl} model). With this, one can aggregate the original problem in different ways to gain insight on the problem structure or to use different types of solution capabilities available (e.g., those that exploit or do not exploit structures).

\subsection{{\tt MadNLP.jl}}

{\tt MadNLP.jl} is a nonlinear programming (NLP) solver implemented in the {\tt Julia} programming language. {\tt MadNLP.jl} implements a state-of-the-art, filter line-search method (as that implemented in {\tt Ipopt}). {\tt MadNLP.jl} has the ultimate goal of facilitating the implementation of algorithmic and linear algebra strategies that can enhance robustness and efficiency. Such flexibility allows for the exploitation of diverse structures present in NLPs (e.g., sparsity patterns or high-level graph structures). Here, we provide a high-level introduction of the interior-point method and of the Schur complement decomposition method implemented in {\tt MadNLP.jl}.

\subsubsection{Interior-Point Method}

In an interior-point method, one attempts to find the solution of \eqref{eq:optimization_model} by solving a sequence of smooth {\it barrier subproblems}. The barrier subproblems are formulated by replacing the inequality constraints by log-barrier terms as:
\begin{subequations}\label{eqn:barrier}
  \begin{align}
    \min_{\bx} &\;\;   f(\bx) - \sum_{i=1}^n \mu \log \bx_{(i)}\\
    \textrm{s.t.} &\; c(\bx) = 0,\quad (\blambda),
  \end{align} 
\end{subequations}
where $\bx_{(i)}$ denotes the $i$-th component of $\bx$. The barrier subproblems are solved using a decreasing sequence of values for the the penalty parameter $\mu>0$. For convenience, we let $\bX:=\diag(\bx)$, $\bZ:=\diag(\bz)$, and $\bone:=[1;1;\cdots;1]$.
\\

To solve \eqref{eqn:barrier}, we compute search (Newton-like) steps by solving a linearized version of the first-order optimality conditions to \eqref{eqn:barrier}. To see this, we write the optimality conditions for \eqref{eqn:barrier}:
\begin{subequations}\label{eqn:barrier-kkt}
  \begin{align}
    \nabla_{\bx} f(\bx) - \mu\bX^{-1} \bone +  \nabla_{\bx}c(\bx)^\top \blambda &= 0\\
    c(\bx)& = 0,
  \end{align}
\end{subequations}
or equivalently, by introducing $\bz = \mu\bX^{-1} \bone$,
\begin{subequations}\label{eqn:barrier-kkt}
  \begin{align}
    \nabla_{\bx} f(\bx) - \bz +  \nabla_{\bx}c(\bx)^\top \blambda &= 0\\
    c(\bx)& = 0\\
    \bX\bz -\mu\bone &= 0.
  \end{align}
\end{subequations}
The Newton step at iteration $k$ is given by:
\begin{align}\label{eqn:barrier-newton}
  \begin{bmatrix}
    \bW^{(k)} & {\bA^{(k)}}^\top & -I\\
    \bA^{(k)}&&\\
    \bZ^{(k)} && \bX^{(k)}\\
  \end{bmatrix}
  \begin{bmatrix}
    \bd^{(k)}_x\\
    \bd^{(k)}_\lambda\\
    \bd^{(k)}_z
  \end{bmatrix}
  =
  \begin{bmatrix}
    \nabla_{\bx} f(\bx^{(k)})   -\bz^{(k)} + {\bA^{(k)}}^\top \blambda^{(k)}\\
    c(\bx^{(k)})\\
    \bX^{(k)} \bZ^{(k)} \bone - \mu \bone
  \end{bmatrix},
\end{align}
where $\bW^{(k)}:=\nabla^2_{\bx\bx} \mathcal{L}(\bx^{(k)},\blambda^{(k)},\bz^{(k)})$, $\bA^{(k)} := \nabla_{\bx}c(\bx^{(k)})$. The right-hand side of \eqref{eqn:barrier-newton} comes from the linearization of \eqref{eqn:barrier-kkt} and the right-hand-side is the evaluation of the optimality system \eqref{eqn:barrier-kkt} at the current iteration. By eliminating the third block-row, we obtain
\begin{align*}
  \begin{bmatrix}
    \bW^{(k)} + \bSigma^{(k)} &{\bA^{(k)}}^\top\\
    \bA^{(k)}
  \end{bmatrix}
  \begin{bmatrix}
    \bd^{(k)}_x\\
    \bd^{(k)}_\lambda
  \end{bmatrix}=
  \begin{bmatrix}
    \nabla_{\bx}f(\bx^{(k)}) - \mu {\bX^{(k)}}^{-1}\bone  + {\bA^{(k)}}^\top \blambda^{(k)}\\
    c(\bx^{(k)})
  \end{bmatrix}
\end{align*}
where $\bSigma^{(k)}:={\bX^{(k)}}^{-1} \bZ^{(k)}$. The solution $(\bd^{(k)}_x,\bd^{(k)}_\lambda,\bd^{(k)}_z)$ of \eqref{eqn:barrier-newton} is used to compute the Newton step. Then, a line-search algorithm is used to determine the steplength $\alpha$ and the iterate:
\begin{align*}
  \bx^{(k+1)} &= \bx^{(k)}+\alpha\cdot \bd^{(k)}_x\\
  \blambda^{(k+1)} &= \blambda^{(k)}+\alpha \cdot \bd^{(k)}_\lambda\\
 \bz^{(k+1)} &= \bz^{(k)}+\alpha \cdot \bd^{(k)}_z.
\end{align*}
The steplength $\alpha$ is chosen so as to improve feasibility or the objective function, see  \cite{wachter2006}. 

\subsubsection{Schur Decomposition}

Graph modeling greatly facilitates the implementation of decomposition strategies; among these, the Schur complement decomposition approach has shown to be flexible and robust.  The details of this approach can be found in \cite{chiang2014,jalving2021,kang2014}. In the context of the graph-structured NLP \eqref{eq:graph_model}, the KKT conditions can be written as:
\begin{subequations}\label{eq:KKT_system}
\begin{align} 
\nabla_{\bx_n} f_n(\bx_n) + \mu \bX^{-1}_n\bone + \nabla_{\bx_n}c_n^I(\bx_n)\blambda_{I,n}  +  \sum_{e \in \mathcal{E}(n)} \nabla_{\bx_{n }} c_e^L(\{\bx_{n'}\}_{n' \in \mathcal{N}(e)})^\top \blambda_{L,e} &= 0,\quad n\in\mathcal{N}\\
c_n^I(\bx_n) &= 0,\quad n \in \mathcal{N} \\
c_e^L(\{\bx_{n}\}_{n \in \mathcal{N}(e)}) &= 0, \quad e \in \mathcal{E},
\end{align} 
\end{subequations}
where $\bX_n:=\diag(\bx_n)$, $\mathcal{E}(n)$ is the set of edges one of whose support is $n$, and $\blambda_{I,n}$ and $\blambda_{L,e}$ are the dual variables corresponding to the internal and the linking constraints respectively. The corresponding KKT system can be formulated as:

\begin{align}\label{eq:BBD_system}
\left[\begin{array}{cccc|c}
\boldsymbol{K}^{(k)}_{1}& &&& {\boldsymbol{B}^{(k)}_{1}}^\top\\
&\boldsymbol{K}^{(k)}_{2}& && {\boldsymbol{B}^{(k)}_2}^\top\\
&&\ddots&&\vdots\\
&&&\boldsymbol{K}^{(k)}_{n}&{\boldsymbol{B}^{(k)}_n}^\top\\\hline
\boldsymbol{B}^{(k)}_1&\boldsymbol{B}^{(k)}_2&\hdots&\boldsymbol{B}^{(k)}_n& \\
\end{array}\right]
\left[\begin{array}{c} \boldsymbol{d}^{(k)}_{1}\\ \boldsymbol{d}^{(k)}_{2}\\ \vdots \\ \boldsymbol{d}^{(k)}_{n}\\ \hline \bd^{(k)}_L \end{array}\right]=
-\left[\begin{array}{c} \boldsymbol{r}^{(k)}_{1}\\ \boldsymbol{r}^{(k)}_{2}\\  \vdots \\ \boldsymbol{r}^{(k)}_{n}\\ \hline \ \boldsymbol{r}^{(k)}_L\ \end{array}\right].
\end{align}
\noindent Here, $\boldsymbol{d}^{(k)}_{n} := ( \bd^{(k)}_{x,n}, \bd^{(k)}_{I,n})$ are the primal-dual steps for all $n \in \mathcal{N}$ where $\bd^{(k)}_{x,n}$ are the primal steps and $\bd^{(k)}_{I,n}$ are the dual steps associated with $\blambda^I_{n,k}$; $\bd^{(k)}_L$ are the dual steps associated with $\blambda^{(k)}_{L,e}$ for $e \in \mathcal{E}$ (the linking constraints); and $\boldsymbol{K}^{(k)}_{n}$, $\boldsymbol{B}^{(k)}_{n}$, $\boldsymbol{r}^{(k)}_{n}$, and $\boldsymbol{r}^{(k)}_L$ are defined as:
\begin{align*}
  \boldsymbol{K}^{(k)}_{n} &:= \left[\begin{array}{cc} \boldsymbol{W}^{(k)}_{n} + \bSigma^{(k)}_{n} & {\boldsymbol{J}^{(k)}_n}^\top \\ \boldsymbol{J}^{(k)}_{n} & \boldsymbol{0} \end{array} \right],\quad n \in \mathcal{N}\\
  \boldsymbol{B}^{(k)}_{n} &:= \left[ \begin{array}{cc} \boldsymbol{Q}^{(k)}_{n} & \boldsymbol{0} \end{array} \right],\quad n \in \mathcal{N}\\
  \boldsymbol{r}^{(k)}_{n} &:= \left[ \begin{array}{c} \nabla_{\bx_{n}} f(\bx^{(k)}_{n}) \\ c_n^I(\bx^{(k)}_{n})  \end{array} \right],\quad n \in \mathcal{N} \\
  \boldsymbol{r}^{(k)}_L &:= \{c_{e}^L(\{\bx^{(k)}_{n}\}_{n \in \mathcal{N}(e)})\}_{e \in \mathcal{E}} .
\end{align*}
\noindent Here, $\boldsymbol{J}^{(k)}_{n} := \nabla_{\bx_n} c_n^I(\bx^{(k)}_{n})$ is the Jacobian of the internal constraints, $\boldsymbol{W}^{(k)}_{n}$ is the Hessian of the Langrangian of \eqref{eq:graph_model}, $\bSigma^{(k)}_n:={\bX^{(k)}_{n}}^{-1} \bZ^{(k)}_{n}$, and $\boldsymbol{Q}^{(k)}_{n} := \nabla_{\bx_n} \{c_e^L (\{\bx^{(k)}_{n'}\}_{n' \in \mathcal{E}(n)})\}_{e \in \mathcal{E}}$ is the Jacobian of the linking constraints. The block-bordered structure in \eqref{eq:BBD_system} naturally results from the graph structure and  facilitates the implementation of the Schur decomposition approach. We also note that the matrices $\boldsymbol{K}^{(k)}_{n}$ and $\boldsymbol{B}^{(k)}_{n}$ for $n \in \mathcal{N}$ may have varying sizes and that such sizes depend on the model construction and on the number of linking constraints. We also highlight that the diagonal blocks in the linear system can embed additional structures (that might arise from hierarchical graph structures). 
\\

Having the block-bordered structure, the decomposition approach proceeds by forming the Schur complement: 
\begin{subequations}\label{Schur_complement}
\begin{align}
    \boldsymbol{S}^{(k)} &= -\sum_{n\in \mathcal{N}} {\boldsymbol{B}^{(k)}_{n}}^\top {\boldsymbol{K}^{(k)}_{n}} ^{-1} \boldsymbol{B}^{(k)}_{n} \label{Schur_complement-1}\\
    \boldsymbol{S}^{(k)} \bd^{(k)}_L &= \sum_{n\in \mathcal{N}} {\boldsymbol{B}^{(k)}_{n}}^\top {\boldsymbol{K}^{(k)}_{n}}^{-1} \boldsymbol{r}^{(k)}_{n} - \boldsymbol{r}^{(k)}_L\label{Schur_complement-2}\\
    \boldsymbol{K}^{(k)}_{n} \bd^{(k)}_{n} &= \boldsymbol{B}^{(k)}_{n} \bd^{(k)}_L - \boldsymbol{r}^{(k)}_{n}, n \in \mathcal{N}.\label{Schur_complement-3}
\end{align}
\end{subequations}
The computation of ${\boldsymbol{B}^{(k)}_{n}}^\top {\boldsymbol{K}^{(k)}_{n}} ^{-1} \boldsymbol{B}^{(k)}_{n} $ can be parallelized since each block operation is  independent. Once the blocks ${\boldsymbol{B}^{(k)}_{n}}^\top {\boldsymbol{K}^{(k)}_{n}} ^{-1} \boldsymbol{B}^{(k)}_{n} $ are computed in parallel, they are assembled into matrix $\boldsymbol{S}^{(k)}$, and the dense linear system \eqref{Schur_complement-2} is solved to obtain the dual step direction for the link constraints. From this step, the block steps are computed in parallel via \eqref{Schur_complement-3}. Parallelization can significantly speed up the step computation, especially when the $\boldsymbol{K}^{(k)}_{n}$ blocks are large and the Schur complement $\boldsymbol{S}^{(k)}$ is small.

\subsection{Model-Solver Interface}

{\tt MadNLP.jl} is designed to interface with OptiGraph object models created in {\tt Plasmo.jl}. A couple of benefits come from interfacing these capabilities: (i) {\tt Plasmo.jl} creates a graph-based model where the structural information of the problem is stored, and such information can then be automatically exploited with {\tt MadNLP.jl} and solved with specialized algorithms, and (ii) the modular structure within the OptiGraph model allows parallel function and derivative evaluations. When used in conjuction with {\tt Plasmo.jl}, {\tt MadNLP.jl} detects the partitions as defined by {\tt Plasmo.jl}. 
\\

If the modeler uses a more general (unstructured) modeling language (such as {\tt JuMP.jl}), a custom partition must be manually set with a vector directing which nodes belong on which partition. Alternatively, if no custom partition is set, {\tt MadNLP.jl} provides capabilities to partition the problem internally using {\tt Metis}. This latter method may be especially useful if the user does not have advanced knowledge of the problem structure or thinks that there may be additional (non-obvious) structures to exploit. 

%%%%%%%%%%%%%%%%%%%%%%%%%%%%%%%%%%%%%%%%%%

\section{Illustrative Example}

In this section, we highlight differences that arise when modeling problems using a general algebraic package such as {\tt JuMP.jl} and the graph-based package {\tt Plasmo.jl}. We also illustrate how graph-based modeling allows users to formulate and partition a problem and note how this can be useful in solving a problem. 

\subsection{Model}

We consider an example that arises in the optimal tuning of a proportional integral derivative (PID) controller. The problem is a stochastic program designed to minimize the error across multiple operational scenarios by choosing tuning parameters. The problem is formulated as:
\begin{subequations}\label{eq:grid}
\begin{align} 
\min_{\{x_s, u_x \}_{s \in S}, K_c, \tau_I, \tau_D}&\;\;   \frac{1}{|S|} \sum_{s\in S} \int_0^{t_f}\left(100(x_{sp,s}-x_s(t))^2 + 0.01 u_s(t)^2\right) \label{eqn:PID_objective_function}\\
\textrm{s.t.} &\;  \frac{1}{\tau}  \frac{dx_s(t)}{dt} + x_s(t) = K u_s(t) + K_d d_s \label{eq:controller_gain}\\
&\; u_s(t) = K_c \left(x_{sp,s}(t) - x_s(t)\right) + \tau_I \int_0^t\left(x_{sp,s}(t) - x_s(t)\right)dt + \tau_D \frac{dx_s(t)}{dt}\label{eq:PID_constraint}\\
&\; x_s(t) = x_0\label{eq:setpoint_constraint}\\
&\; -10 \leq K_c \leq 10 \label{eq:limits_constraint_start}\\
&\; -100 \leq \tau_I \leq 100\\
&\; -100 \leq \tau_D \leq 100\\
&\; -2.5 \leq x_s(t) \leq 2.5\\
&\; -2.0 \leq u_s(t) \leq 2.0 \label{eq:limits_constraint_end}.
\end{align} 
\end{subequations}
\noindent Here, $t$ is the time for $t \in [0,t_f]$ where $t_f = 10$, and $s$ is the scenario index for $s\in S$, where $S$ is a set of five scenarios. The state variable is represented by $x_s(t)$ for $s \in S$, and the set-point is represented by $x_{sp,s}$ for $s \in S$. Equation \eqref{eqn:PID_objective_function} is the objective function penalizing the deviation from the set-point and the magnitude of the control actions. Equation \eqref{eq:controller_gain} is the constraint representing a first-order linear dynamical system with a disturbance and without time delay. Here,   $\tau$ is the time constant, $K$ is the process gain, $K_d$ is the disturbance gain, and $d_s$ is the disturbance for $s \in S$. Equation \eqref{eq:PID_constraint} is the PID controller formulation, where the tuning parameters $K_c$, $\tau_I$, and $\tau_D$ relate to the proportional, integral, and derivative terms. These tuning parameters are the design variables and will thus have the same value across scenarios. For simplicity, \eqref{eq:setpoint_constraint} requires all scenarios to have the same starting point ($x_0$). Equations \eqref{eq:limits_constraint_start}-\eqref{eq:limits_constraint_end} give upper and lower limits for all decision variables. 

\subsection{Implementation}

In implementing the problem as a {\tt JuMP.jl} model object or a  {\tt Plasmo.jl} OptiGraph object, we discretize the dynamical system using 100 evenly-spaced time points . Scenario and operating parameters are shown in the code snippet of Figure \ref{fig:code_snippet_PID_IC}. 
\\

The formulation as a {\tt JuMP.jl} object is given in the code snippet in Figure \ref{fig:code_snippet_PID_JuMP}.  Comments within the snippet explain different aspects of the implementation. We highlight that this model uses indices for time and scenario to define process variables or inputs; moreover, \textit{this formulation has no predefined structure}. 
\\

The corresponding formulation in {\tt Plasmo.jl} is shown in the code snippet of Figure \ref{fig:code_snippet_PID_Plasmo} and contains an inherent structure by defining the problem with nodes and edges. In this construction, nodes represent time points within each scenario. A master node is added to contain the design variables (i.e. the controlled tuning parameters) that are then linked across scenarios. In this formulation, dummy variables are introduced (lines \ref{line:dummy_start} - \ref{line:dummy_end}) to avoid nonlinearities in the linking constraints (this is a lifting procedure). There are several ways in which one can formulate the tuning problem in {\tt Plasmo.jl}; for instance, the dummy variables could have been avoided by placing the nonlinear constraints on a single node (placing a full scenario in each node). However, the more variables and constraints contained in a single node, the fewer options that exist for partitioning the problem into a  hierarchical subgraphs. Here, we chose to model this problem with a large number of nodes to allow for additional partitioning. This is an important consideration when formulating a model using a graph structure. 

\begin{figure}[!htp]
\centering
\begin{scriptsize}
\lstset{language=Julia, breaklines = true, xleftmargin=\parindent}
\begin{lstlisting}[escapeinside={(*}{*)}, escapechar=|]
     using Plasmo, JuMP
     
     # sets
     NS= 5      # number of scenarios
     N=100;     # number of timesteps
     Tf=10;     # final time
     h=Tf/N;    # time step
     T=1:N;     # set of times
     Tm=1:N-1;  # set of times minus one
     
     # set time vector
     time=zeros(N);
     for t=1:N
       time[t] = h*(t-1);
     end
     
     # Define parameters used within each scenario
     K    = 1.0;      # Gain
     x0   = 0.0;      # Starting Point
     Kd   = 0.5;      # Disturbance Gain
     tau  = 1.0;      # Time Constant
     d    = fill(-1.0,5);      # Disturbance
     xsp  = [-2.0,-1.5,-0.5,0.5,1.0] # Set Point
\end{lstlisting}
\end{scriptsize}
\caption{Code Snippet setting parameters for the PID controller tuning problem}
\label{fig:code_snippet_PID_IC}
\end{figure}

\begin{figure}[!htp]
\centering
\begin{scriptsize}
\lstset{language=Julia, breaklines = true, xleftmargin=\parindent}
\begin{lstlisting}[escapeinside={(*}{*)}, escapechar=|]
# define JuMP model
m = Model()

# define model variables
@variable(m, -2.5 <= x[S, T] <= 2.5) # process variable
@variable(m, -2.0 <= u[S, T] <= 2.0) # process input
@variable(m, int[S, T]) # discretized integral at time t
@variable(m, cost[S,T]) # cost to objective function

# define tuning parameters (first stage variables)
@variable(m, -10 <= Kc <= 10)
@variable(m, -100 <= tauI <= 100)
@variable(m, -100 <= tauD <= 100)
    
# constrain model to follow first order linear system
@constraint(m, [s in S, t in Tm], (1/tau)*(x[s,t+1]-x[s,t])/h + x[s,t+1] == K*u[s,t+1] + Kd*d[s]); 

# define controller operation
@constraint(m, [s in S, t in Tm], u[s,t+1] == Kc*(xsp[s] - x[s,t]) 
															 + tauI*int[s,t+1] + tauD*(x[s,t+1] - x[s,t])/h);
# discretize integral
@constraint(m, [s in S, t in Tm], (int[s,t+1] - int[s,t])/h == xsp[s] - x[s,t+1]);

# set initial conditions
@constraint(m, [s in S], x[s,1] == x0);
@constraint(m, [s in S], int[s,1] == 0);
# define cost for objective function
@constraint(m, eqcost[s in S, t in T], cost[s,t] == 100*(xsp[s]-x[s,t])^2 + 0.01*u[s,t]^2);

#objective function
@objective(m, Min, (1/(NS)) * sum(cost[s,t] for s in S, t in T))
\end{lstlisting}
\end{scriptsize}
\caption{Code Snippet defining a PID controller problem formulated in {\tt JuMP.jl}}
\label{fig:code_snippet_PID_JuMP}
\end{figure}

\begin{figure}[!htp]
\centering
\begin{scriptsize}
\lstset{language=Julia, breaklines = true, xleftmargin=\parindent}
\begin{lstlisting}[escapeinside={(*}{*)}, escapechar=|]
     # define a function to create a graph for each scenario
     function get_scenario_model(s)
       graph = OptiGraph() # create OptiGraph
     
       @optinode(graph,n[T]) # add N nodes, n[1:100], to the optigraph
     
       # add decision variables to each node
       for node in n
     
         # define variables on node n
         @variable(node,-2.5<= x <=2.5)
         @variable(node,-2.0<= u <=2.0)
         @variable(node, int )
         @variable(node, Kc)
         @variable(node, tauI)
         @variable(node, tauD)
     
         # define dummy variables to avoid nonlinearity in @linkconstraint
         @variable(node, Kcx) |\label{line:dummy_start}|
         @variable(node, tauIint)
         @variable(node, tauDx)
         
         # define dummy variable constraints
         @constraint(node, Kcx == Kc * x)
         @constraint(node, tauIint == tauI * int)
         @constraint(node, tauDx  == tauD * x) |\label{line:dummy_end}|
     
         # objective function
         @objective(node, Min, (1/(NS))*(100*(xsp[s]-x)^2 + 0.01*u^2))
       end
     
       # set initial conditions
       @constraint(n[1], eqinix,  n[1][:int] == 0)
       @constraint(n[1], eqinit,  n[1][:x] == x0)
     
       # constraints
       @linkconstraint(graph, [t=Tm], (1/tau)*(n[t+1][:x]-n[t][:x])/h + n[t+1][:x]
                                                                                  == K*n[t+1][:u] + Kd*d[s]);

       @linkconstraint(graph, [t=Tm], n[t+1][:u] == n[t][:Kc]*xsp[s] - n[t][:Kcx] 
                                                  + n[t+1][:tauIint] + n[t+1][:tauDx]/h - n[t][:tauDx]/h);

       @linkconstraint(graph, [t=Tm], (n[t+1][:int]-n[t][:int])/h == (xsp[s]-n[t+1][:x]));
     
       #link stage one decision variables
       @linkconstraint(graph, [t=Tm], n[t][:Kc] == n[t+1][:Kc])
       @linkconstraint(graph, [t=Tm], n[t][:tauI] == n[t+1][:tauI])
       @linkconstraint(graph, [t=Tm], n[t][:tauD] == n[t+1][:tauD])
     
       return graph, n
     end

     PID=OptiGraph() # create initial high-level OptiGraph
     
     @optinode(PID,master) # create node called "master" for first stage variables
     
     # define tuning parameters (first stage variables)
     @variable(master, -10<= Kc <=10)
     @variable(master,-100<= tauI <=100)
     @variable(master,-100<= tauD <=100)
     
     for s in 1:NS           
       graph,nodes = get_scenario_model(s) # get OptiGraph and OptiNode for each Scenario

       add_subgraph!(PID,graph) # add scenario subgraph to high-level OptiGraph |\label{line:add_subgraphs}|
     
       # link constraints from subgraph to master node
       @linkconstraint(PID, nodes[1][:Kc]==master[:Kc])
       @linkconstraint(PID, nodes[1][:tauI]==master[:tauI])
       @linkconstraint(PID, nodes[1][:tauD]==master[:tauD])    
     end
\end{lstlisting}
\end{scriptsize}
\caption{Code Snippet defining a PID controller problem formulated in {\tt Plasmo.jl}}
\label{fig:code_snippet_PID_Plasmo}
\end{figure}

\subsection{Partitioning}

\begin{figure}
     \centering
     \includegraphics[scale=.6]{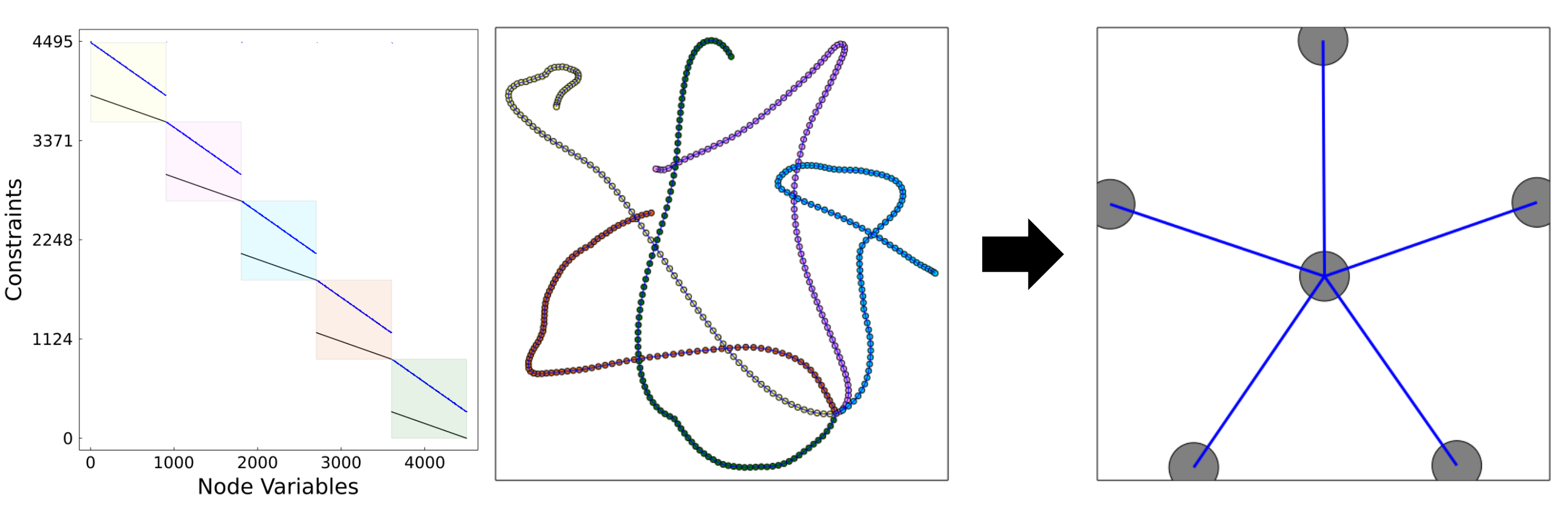}
     \caption{Visualization of the PID controller problem partitioned by scenario. Adjacency matrix (left) is shown with the graph-based model partitioned by scenario with colors representing each subgraph (middle) and the aggregated form of the graph (right)}
     \label{fig:Scen_Figure}
 \end{figure}

The problem in Figure \ref{fig:code_snippet_PID_Plasmo} is partitioned based on scenarios, because each individual scenario is placed on a separate subgraph (Line \ref{line:add_subgraphs}). The high-level OptiGraph consists of a master node connected to five separate subgraphs. Using the {\tt aggregate} function, each subgraph can be combined into a single node (creating a new OptiGraph). Visualizations of the scenario-partitioned model in its aggregated and non-aggregated form, along with their accompanying adjacency matrices, are shown in Figure \ref{fig:Scen_Figure}. These were created using {\tt Plasmo.jl} visualization capabilities  \cite{jalving2021}. These visualizations reveal the tree structure seen in typical stochastic optimization formulations and highlight how the modeler can use visualization capabilities to navigate and understand the problem structure.
 
\begin{figure}
     \centering
     \includegraphics[scale=.6]{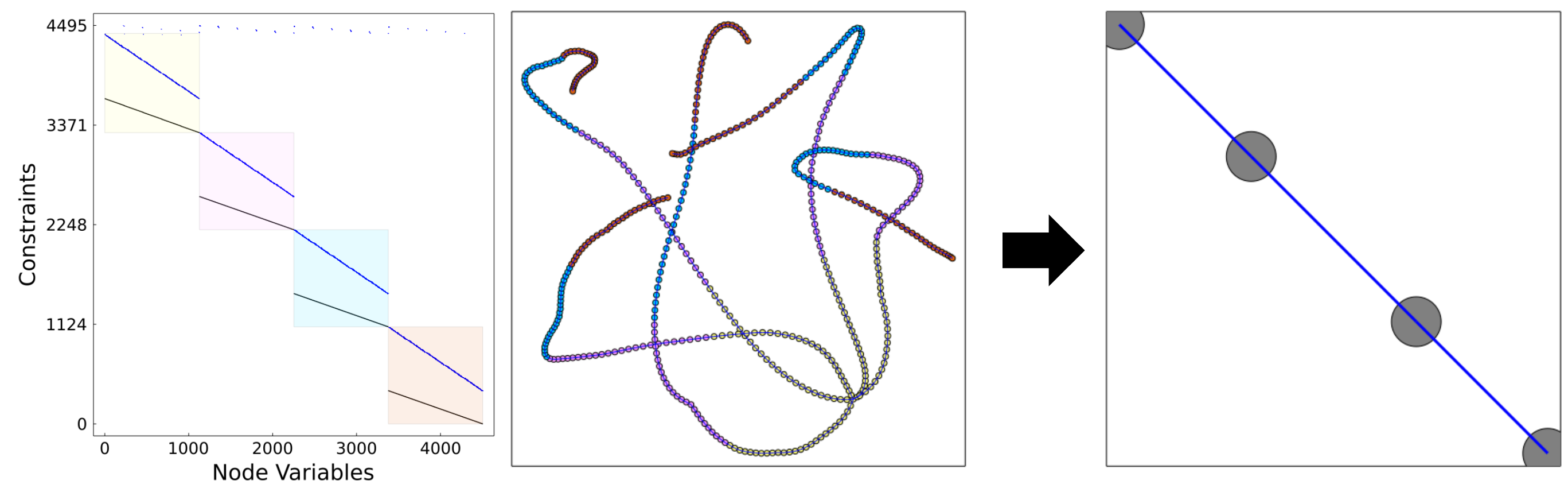}
     \caption{Visualization of the PID controller problem partitioned by time. Adjacency matrix (left) is shown with the graph-based model partitioned by scenario with colors representing each subgraph (middle) and the aggregated form of the graph (right)}
     \label{fig:Time_Figure}
 \end{figure}

With the initial OptiGraph created, the problem can be partitioned in other ways. Desired partitioning schemes will depend on the model structure and intended decomposition method. For the PID controller, we show that the problem can also be partitioned in the time domain (rather than the scenario domain). The code snippet in Figure \ref{fig:code_snippet_PID_partition} shows how the formulation in Figure \ref{fig:code_snippet_PID_Plasmo} can be repartitioned over time. Each new partition contains 25 time points from each scenario. We first create a reference map (line \ref{line:refmap}) and a node membership vector (line \ref{line:nmv}). The node membership vector contains an integer value corresponding to each node of the reference map that directs {\tt Plasmo.jl} into which partition to place each node. When the membership vector is filled (lines \ref{line:nmv_fill_start} - \ref{line:nmv_fill_end}), a new partition can be created (line \ref{line:make_partition}) and applied to the original OptiGraph (line \ref{line:make_subgraphs}). The resulting OptiGraph is partitioned in time and is visualized in Figure \ref{fig:Time_Figure}, along with the new adjacency matrix and its corresponding aggregated structure. This structure reveals the typical linear tree structure of a dynamic optimization problem. 

\begin{figure}
\centering
\begin{scriptsize}
\lstset{language=Julia, breaklines = true, xleftmargin=\parindent}
\begin{lstlisting}[escapeinside={(*}{*)}, escapechar=|]
     # create reference map for partition
     hypergraph, refmap = gethypergraph(PID) |\label{line:refmap}|
     
     # define a node membership vector that assigns each node index to a partition
     node_membership_vector = Array{Int64,1}(undef,length(all_nodes(PID))) |\label{line:nmv}|
     
     # fill node membership vector; value of vector entry corresponds to partition number
     node_membership_vector[[1:26; 102:126; 202:226; 302:326; 402:426]]   .= 1 |\label{line:nmv_fill_start}|
     node_membership_vector[[27:51; 127:151; 227:251; 327:351; 427:451]]  .= 2
     node_membership_vector[[52:76; 152:176; 252:276; 352:376; 452:476]]  .= 3
     node_membership_vector[[77:101; 177:201; 277:301; 377:401; 477:501]] .= 4 |\label{line:nmv_fill_end}|
     
     # create partition
     PID_partition = Partition(node_membership_vector,refmap) |\label{line:make_partition}|
     
     # repartition subgraphs according to the partition
     make_subgraphs!(PID, PID_partition) |\label{line:make_subgraphs}|
\end{lstlisting}
\end{scriptsize}
\caption{Code Snippet for creating manual partition of an OptiGraph}
\label{fig:code_snippet_PID_partition}
\end{figure}

Partitioning a problem can result in quite different structures. Figures \ref{fig:Scen_Figure} and \ref{fig:Time_Figure} contain the same problem but present distinct structures. These structures may be best exploited by using certain decomposition schemes. For example, the aggregated structure shown in Figure \ref{fig:Scen_Figure} can be solved efficiently with Schur decomposition because it is a bilevel tree. Alternatively, the aggregated structure seen in Figure \ref{fig:Time_Figure} may best be exploited with Schwarz decomposition because subproblems could have a degree of overlap without encompassing the whole graph. These modeling and partitioning capabilities of graph-structured optimization can thus be useful for analyzing and solving problems, as evidenced above.

%%%%%%%%%%%%%%%%%%%%%%%%%%%%%%%%%%%%%%%%%%
\section{Large-Scale Case Study}

\begin{figure}[!htp]
  \centering
  \includegraphics[scale=.65]{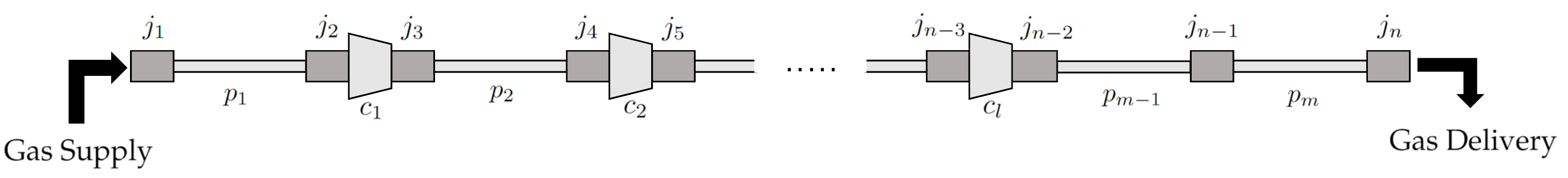}
  \caption{Visualization of a natural gas pipeline containing $l$ compressors, $m$ pipelines, and $n$ junctions. Junctions are placed at connections between two pipelines and connections between a pipeline and compressor. $c_i$ represents the $i$th compressor, $p_i$ the $i$th pipeline, and $j_i$ the $i$th junction. Gas is supplied at the first junction and delivered at the $n$th junction.}
  \label{fig:Gas_Pipeline}
\end{figure}

To demonstrate the scalability of {\tt Plasmo.jl} and {\tt MadNLP.jl}, we implemented a large-scale natural gas network problem. This is a two-stage stochastic optimization problem for a system that comprises pipelines, compressors, and junctions (Figure \ref{fig:Gas_Pipeline}). The first-stage variables are given by the  compressor power policy, while the second-stage (recourse) variables are the states of the pipelines (e.g., pressures and flows) in different scenarios. Each scenario within the problem contained a different gas demand profile at the end of the network. 

\subsection{Model}

To construct this problem, we used the constraints defined by \cite{jalving2021, jalving2018} and the code for a deterministic problem given by \cite{jalving2021} and modified the code to  capture uncertainty. There are three object types that lead to modeling constraints: junctions, compressors, and pipelines. There are also additional constraints created by links between these object types. The constraints and objective function are given in subsections \ref{sec:junction_constraints} - \ref{sec:obj_function} below. The formulations presented below follow closely those defined in \cite{jalving2021}. However, we introduce a few minor changes to form a stochastic model and to more closely follow the code used to construct this model.  
\\

The following sets are employed within the given formulation: $\mathcal{L}_c$ represents the set of all compressors, $\mathcal{L}_p$ represents the set of all pipelines, and $\mathcal{L} := \mathcal{L}_c \cup \mathcal{L}_p$ represents the union of the set of compressors and set of pipelines. We also use $\mathcal{J} := \{1, ..., N_j \}$ for the set of all junctions in the network, where a junction is placed at any connection between pipelines and between any connection of a pipeline with a compressor. Junction $1$ is where gas is supplied to the network and junction $N_j$ is the final junction of the pipeline and the location of gas delivery. We further define the set $\mathcal{J}' := \mathcal{J} \backslash \{1, N_j\}$. Pipelines were also discretized into a set of equally spaced points, $\mathcal{X} := \{1, ..., N_x\}$, and the operating time was discretized into a set of time points, $\mathcal{T} := \{1, ..., N_t\}$. Here, $N_x$ and $N_t$ are the final spacial and time points respectively. We also define the sets $\mathcal{X}' := \mathcal{X} \backslash N_x$ and $\mathcal{T}' := \mathcal{T} \backslash 1$. Lastly, we use $\mathcal{S}$ to represent the set of all scenarios.\\

\subsubsection{Junction Constraints}\label{sec:junction_constraints}

The junction constraints are given in Equation \eqref{eqn:junction_constraints}. $\theta_{j,t,s}$ is the pressure for $j \in \mathcal{J}$, $t \in \mathcal{T}$, and $s \in \mathcal{S}$, with $\ubar{\theta}_j$ as the lower bound and $\bar{\theta}_j$ as the upper bound for $j \in \mathcal{J}$. $\bar{F}_{t,s}$ is the amount of supplied gas at the first junction for $t \in \mathcal{T}$ and  $s \in \mathcal{S}$ and bounded by a maximum supply amount, $\bar{\phi}_{supply}$. $F_{t,s}$ is the amount of gas being delivered at the end of the network for $t \in \mathcal{T}$ and $s \in \mathcal{S}$. $\tilde{F}_{t,s}$ is a positive variable that allows for $F_{t,s}$ to exceed the demand ($d_{t,s}$) within a given scenario $s$ at time $t$ (Equation \eqref{eq:slack}). $\tilde{F}_{t,s}$ is also penalized within the objective function (Equation \eqref{eq:objective_function}) to help limit over supplying gas in any scenario.

\begin{subequations}\label{eqn:junction_constraints}
  \begin{align}
    & F_{t,s} \le d_{t,s} + \tilde{F}_{t,s}, \qquad t \in \mathcal{T}, s \in \mathcal{S} \label{eq:slack}\\
    & \ubar{\theta}_j \le \theta_{j,t,s} \le \bar{\theta}_j, \qquad j \in \mathcal{J}, t \in \mathcal{T}, s \in \mathcal{S}\\
    & 0 \le \bar{F}_{t,s} \le \bar{\phi}_{supply}, \qquad t \in \mathcal{T}, s \in \mathcal{S}\\
    & 0 \le F_{t, s} , \qquad t \in \mathcal{T}, s \in \mathcal{S}\\
    & 0 \le \tilde{F}_{t,s}, \qquad t \in \mathcal{T}, s \in \mathcal{S}
  \end{align}
\end{subequations}

\subsubsection{Compressor Constraints}

The compressor constraints are given in Equation \eqref{eqn:compressor_constraints}. Here,  $p^{in}_{\ell, t, s}$ is the suction pressure, $p^{out}_{\ell, t,s}$ is the discharge pressure, $\bar{p}_{\ell, t,s}$ is the boost in pressure between the suction and discharge , $P_{\ell, t,s}$ is the power of the compressor, and $f_{\ell,t,s}$ is the flow passing through the compressor for $\ell \in \mathcal{L}_c$, $t \in \mathcal{T}$, and $s \in \mathcal{S}$. The parameters $c_P$, $T$, and $\gamma$ are the heat capactiy, temperature, and isentropic efficiency, respectively. The parameters $\ubar{\theta}_{in}$, $\ubar{\theta}_{out}$, $\ubar{\theta}_{boost}$, $\bar{\theta}_{in}$, $\bar{\theta}_{out}$, $\bar{\theta}_{boost}$, and $\bar{\psi}_{power}$ represent lower or upper limits on their respective variables. The operation of the compressor is dictated by Equation \eqref{eqn:comp1}. Lastly, we note that the definition of $f^{in}_{\ell,t,s}$ and $f^{out}_{\ell,t,s}$ for $\ell \in \mathcal{L}_c, t \in \mathcal{T}, s \in \mathcal{S}$ in Equation \eqref{eqn:comp_links} is done to simplify the mathematical definition of the linking constraints (see Equation \eqref{eqn:link_constraints}).

\begin{subequations}\label{eqn:compressor_constraints}
  \begin{align}
    & P_{\ell, t,s} = c_P \cdot T \cdot f_{\ell, t,s} \left[ \left(\frac{p^{out}_{\ell,t,s}}{p^{in}_{\ell,t,s}}\right)^{\frac{\gamma - 1}{\gamma}} -1 \right], \qquad \ell \in \mathcal{L}_c, t \in \mathcal{T}, s \in \mathcal{S} \label{eqn:comp1} \\
    & p^{out}_{\ell,t,s} = p^{in}_{\ell,t,s} + \bar{p}_{\ell,t,s}, \qquad \ell \in \mathcal{L}_c, t \in \mathcal{T}, s \in \mathcal{S}\label{eqn:comp2} \\
    & 0 \le P_{\ell, t, s} \le \bar{\psi}_{power}, \qquad \ell \in \mathcal{L}_c, t \in \mathcal{T}, s \in \mathcal{S} \label{eqn:power_limit}\\
    & \ubar{\theta}_{in} \le p^{in}_{\ell,t,s} \le \bar{\theta}_{in}, \qquad \ell \in \mathcal{L}_c, t \in \mathcal{T}, s \in \mathcal{S}\\
    & \ubar{\theta}_{out}\le p^{out}_{\ell,t,s} \le \bar{\theta}_{out},\qquad \ell \in \mathcal{L}_c, t \in \mathcal{T}, s \in \mathcal{S}\\
    & \ubar{\theta}_{boost} \le \bar{p}_{\ell,t,s} \le \bar{\theta}_{boost}, \qquad \ell \in \mathcal{L}_c, t \in \mathcal{T}, s \in \mathcal{S}\\
    & 0 \le f_{\ell,t,s} , \qquad \ell \in \mathcal{L}_c, t \in \mathcal{T}, s \in \mathcal{S}\\
    & f_{\ell,t,s} = f^{in}_{\ell,t,s} = f^{out}_{\ell,t,s}, \qquad \ell \in \mathcal{L}_c, t \in \mathcal{T}, s \in \mathcal{S} \label{eqn:comp_links}
  \end{align}
\end{subequations}

\subsubsection{Pipeline Constraints}

The pipeline constraints are definied in Equation \eqref{eqn:pipeline_constraints}. Here, $p_{\ell, k,t,s}$ and $f_{\ell,k,t,s}$ are the pressure and flows respectively for $\ell \in \mathcal{L}_P$, $k \in \mathcal{X}$, $t \in \mathcal{T}$, and $s \in \mathcal{S}$. We also define the constants $c_{1,\ell}$, $c_{2,\ell}$, and $c_{3,\ell}$ for $\ell \in \mathcal{L}_p$ which relate to the physical properties of the gas and pipelines (see \cite{jalving2021}). $\Delta t$ is the size of the discretized time intervals and $\Delta x_{\ell}$ is the size of hte discretized space interval on pipeline $\ell \in \mathcal{L}_p$. We also introduce the variables $m_{\ell, t,s}$ for $\ell \in \mathcal{L}_p, t \in \mathcal{T}$, and $s \in \mathcal{S}$ as the line pack of a given pipeline. Equations \eqref{eqn:pipe1} and \eqref{eqn:pipe2} relate to the mass and momentum dynamics of the pipeline, Equations \eqref{eqn:ss1} and \eqref{eqn:ss2} require the problem to start at a steady state, Equations \eqref{eqn:linepack1} and \eqref{eqn:linepack2} define the linepack and require the pipelines to refill at the end of the time period to at least their initial linepack, and Equations \eqref{eqn:pipelinks1} - \eqref{eqn:pipelinks2} define flow and pressure terms that are used within the mathematical definition fo the linking constraints (see Equation \eqref{eqn:link_constraints}). 

\begin{subequations}\label{eqn:pipeline_constraints}
  \begin{align}
    & \frac{p_{\ell, t,k, s} - p_{\ell,t-1, k,s}}{\Delta t} + c_{1,\ell} \frac{f_{\ell,  t,k+1, s} - f_{\ell, t,k, s}}{\Delta x_\ell} = 0, \qquad \ell \in \mathcal{L}_p,  t \in \mathcal{T}',k \in \mathcal{X}', s \in \mathcal{S} \label{eqn:pipe1}\\
    & \frac{f_{\ell, t,k, s} - f_{\ell,t-1, k,s}}{\Delta t} = -c_{2,\ell} \frac{p_{\ell, t, k+1,s} - p_{\ell, t, k, s}}{\Delta x_\ell} - c_{3,\ell} \frac{f_{\ell, t,k,s} |f_{\ell,t,k,s}|}{p_{\ell, t,k,s}}, \qquad \ell \in \mathcal{L}_p,  t \in \mathcal{T}',k \in \mathcal{X}', s \in \mathcal{S}\label{eqn:pipe2}\\
    & \frac{f_{\ell, 1, k+1, s} - f_{\ell, 1, k, s}}{\Delta x_\ell} = 0, \qquad \ell \in \mathcal{L}_p, k \in \mathcal{X}', s \in \mathcal{S}\label{eqn:ss1}\\
    & c_{2, \ell} \frac{p_{\ell, 1, k+1,s} - p_{\ell, 1, k, s}}{\Delta x_\ell} + c_{3,\ell} \frac{f_{\ell, 1, k, s} |f_{\ell, 1, k, s}|}{p_{\ell, 1, k, s}} = 0\label{eqn:ss2}\\
    & m_{\ell, t, s} = \frac{1}{c_{1,\ell}} \sum_{k=1}^{N_x - 1} p_{\ell, t, k, s} \Delta x_{\ell}, \qquad \ell \in \mathcal{L}_p, t \in \mathcal{T}, s \in \mathcal{S}\label{eqn:linepack1}\\
    & m_{\ell, N_t,s} \ge m_{\ell, 1, s}, \qquad \ell \in \mathcal{L}_p, s \in \mathcal{S} \label{eqn:linepack2}\\
    & f_{\ell, t, 1, s} = f^{in}_{\ell, t, s}, \qquad \ell \in \mathcal{L}_p, t \in \mathcal{T}, s \in \mathcal{S}\label{eqn:pipelinks1} \\
    & f_{\ell, t, N_x, s} = f^{out}_{\ell, t, s}, \qquad \ell \in \mathcal{L}_p, t \in \mathcal{T}, s \in \mathcal{S} \\
    & p_{\ell, t, 1, s} = p^{in}_{\ell, t, s}, \qquad \ell \in \mathcal{L}_p, t \in \mathcal{T}, s \in \mathcal{S}\\
    & p_{\ell, t, N_x, s} = p^{out}_{\ell, t, s}, \qquad \ell \in \mathcal{L}_p, t \in \mathcal{T}, s \in \mathcal{S} \label{eqn:pipelinks2}
  \end{align}
\end{subequations}

\subsubsection{Linking Constraints}\label{sec:linking_constraints}
Here, we introduce the notation used by \cite{jalving2021} of $\mathcal{L}_{rec(j)}$ as the set of pipelines or compressors that flow into junciton $j$ and $\mathcal{L}_{snd(j)}$ as the set of pipelines or compressors that receive the flow from junction $j$. Further, we define $\theta_{rec(\ell),t,s}$ and $\theta_{snd(\ell),t,s}$ as the receiving and sending junctions respectively for $\ell \in \mathcal{L}$, $t \in \mathcal{T}$, and $s \in \mathcal{S}$. Equations \eqref{eqn:flow_links} - \eqref{eqn:delivery_link} require that the flow going in and out of a junction is conserved and Equations \eqref{eqn:press_link1} and \eqref{eqn:press_link2} require that the pressure is consistent with the objects to which the junction is connected.

\begin{subequations}\label{eqn:link_constraints}
  \begin{align}
    & \sum_{\ell \in \mathcal{L}_{rec(j)}} f^{out}_{\ell, t, s} - \sum_{\ell \in \mathcal{L}_{snd(j)}} f^{in}_{\ell,t,s} = 0, \qquad j \in \mathcal{J}', t \in \mathcal{T}, s \in \mathcal{S} \label{eqn:flow_links}\\
    & \bar{F}_{t,s} - \sum_{\ell \in \mathcal{L}_{snd(1)}} f^{in}_{\ell,t,s} = 0, \qquad  t \in \mathcal{T}, s \in \mathcal{S} \label{eqn:supply_link}\\
    & \sum_{\ell \in \mathcal{L}_{rec(N_j)}} f^{out}_{\ell, t, s} - F_{t,s}, \qquad  t \in \mathcal{T}, s \in \mathcal{S}\label{eqn:delivery_link}\\
    & p^{in}_{\ell, t, s} = \theta_{rec(\ell), t, s}, \qquad \ell \in \mathcal{L}, t \in \mathcal{T}, s \in \mathcal{S} \label{eqn:press_link1}\\
    & p^{out}_{\ell, t, s} = \theta_{snd(\ell), t, s}, \qquad \ell \in \mathcal{L}, t \in \mathcal{T}, s \in \mathcal{S} \label{eqn:press_link2}
  \end{align}
\end{subequations}

\subsubsection{Overall Model}\label{sec:obj_function}

The overall model incorporating the above constraints is given in Equation \eqref{eq:gas_graph_model}. The objective function (Equation \eqref{eq:objective_function}) maximizes the overall profit by taking the minimum of the negative scaled sum of gas delivered plus a scaled compressor power. There is also a penalty term for any amount delivered over a scenario's demand. Positive scalars $\alpha$, $\beta$, and $\kappa$ are weights on their respective terms to give an overall profit. Equations \eqref{eqn:constraint1} - \eqref{eqn:constraint2} are the constraints detailed in sections \ref{sec:junction_constraints} - \ref{sec:linking_constraints} above. We also introduce the variables $\bar{P}_{\ell, t}$ for $\ell \in \mathcal{L}_c$ and $t \in \mathcal{T}$ which are the first stage variables of compressor power. These variables fix the compressor power across all scenarios by requiring that any given compressor operates with the same power at each time point in each scenario (Equation \eqref{eqn:first_stage_variables}).

\begin{subequations}\label{eq:gas_graph_model}
  \begin{align} 
  \min_{\bx} &\;\;   \sum_{\ell \in \mathcal{L}_c, t \in \mathcal{T}, s \in \mathcal{S}} \alpha P_{\ell, t, s} - \sum_{t \in \mathcal{T}, s \in \mathcal{S}} \beta F_{t,s} + \sum_{t \in \mathcal{T}, s \in \mathcal{S}} \kappa \tilde{F}_{t,s}  \label{eq:objective_function}\\
  \textrm{s.t.} &\; \text{Junction Constraints} \quad (\ref{eqn:junction_constraints}) \label{eqn:constraint1}\\
  &\; \text{Compressor Constraints} \quad (\ref{eqn:compressor_constraints})\\
  &\; \text{Pipeline Constraints} \quad (\ref{eqn:pipeline_constraints})\\
  &\; \text{Linking Constraints} \quad (\ref{eqn:link_constraints})\label{eqn:constraint2}\\
  &\; \bar{P}_{\ell, t}  = P_{\ell, t, s}, \qquad \ell \in \mathcal{L}_c, t \in \mathcal{T}, s \in \mathcal{S} \label{eqn:first_stage_variables}
  \end{align} 
\end{subequations}

Finally, note that there are several decision variables within this formulation. In Equation \ref{eq:objective_function}, we introduce the variable $\bx$, which is a vector of all decision variables, including $F_{t,s}$, $\tilde{F}_{t,s}$, $\bar{F}_{t,s}$, $\theta_{j,t,s}$, $p_{\ell, t, s} $, $f_{\ell,t,s}$, $\bar{p}_{\ell,t,s}$, $P_{\ell, t, s}$, $f_{\ell', t,k,s}$, $p_{\ell',t,k,s}$, and $\bar{P}_{t,s}$ for  $\ell \in \mathcal{L}_c$, $\ell' \in \mathcal{L}_p$, $j \in \mathcal{J}$, $t \in \mathcal{T}$, $k \in \mathcal{X}$, and $s \in \mathcal{S}$.

\subsection{Graph Structure and Visualizations}

The gas network formulation forms a graph that places each scenario on its own subgraph within {\tt Plasmo.jl}. These scenario subgraphs were further structured by representing each $\ell \in \mathcal{L}_c$, $\ell' \in \mathcal{L}_p$, and $j \in \mathcal{J}$ as a subgraph. The subgraphs for compressors and junctions were comprised of individual nodes for each $t \in \mathcal{T}$. Pipeline subgraphs were also structured in space, resulting in nodes for each $k \in \mathcal{X}$ in addition to $t \in \mathcal{T}$. In our problem formulation, we used 11 compressors, 13 pipelines, and 25 junctions. We applied discretization to the model to obtain 24 points in time and 10 evenly-spaced points in space for each pipeline. The problem was scaled by changing the number of scenarios within the problem and included as many as 150 scenarios. With each scenario containing 11,376 variables, the largest form of this problem with 150 scenarios contained over 1.7 million variables. 
\\

\begin{figure}[!htp]
     \centering
     \includegraphics[scale=.4]{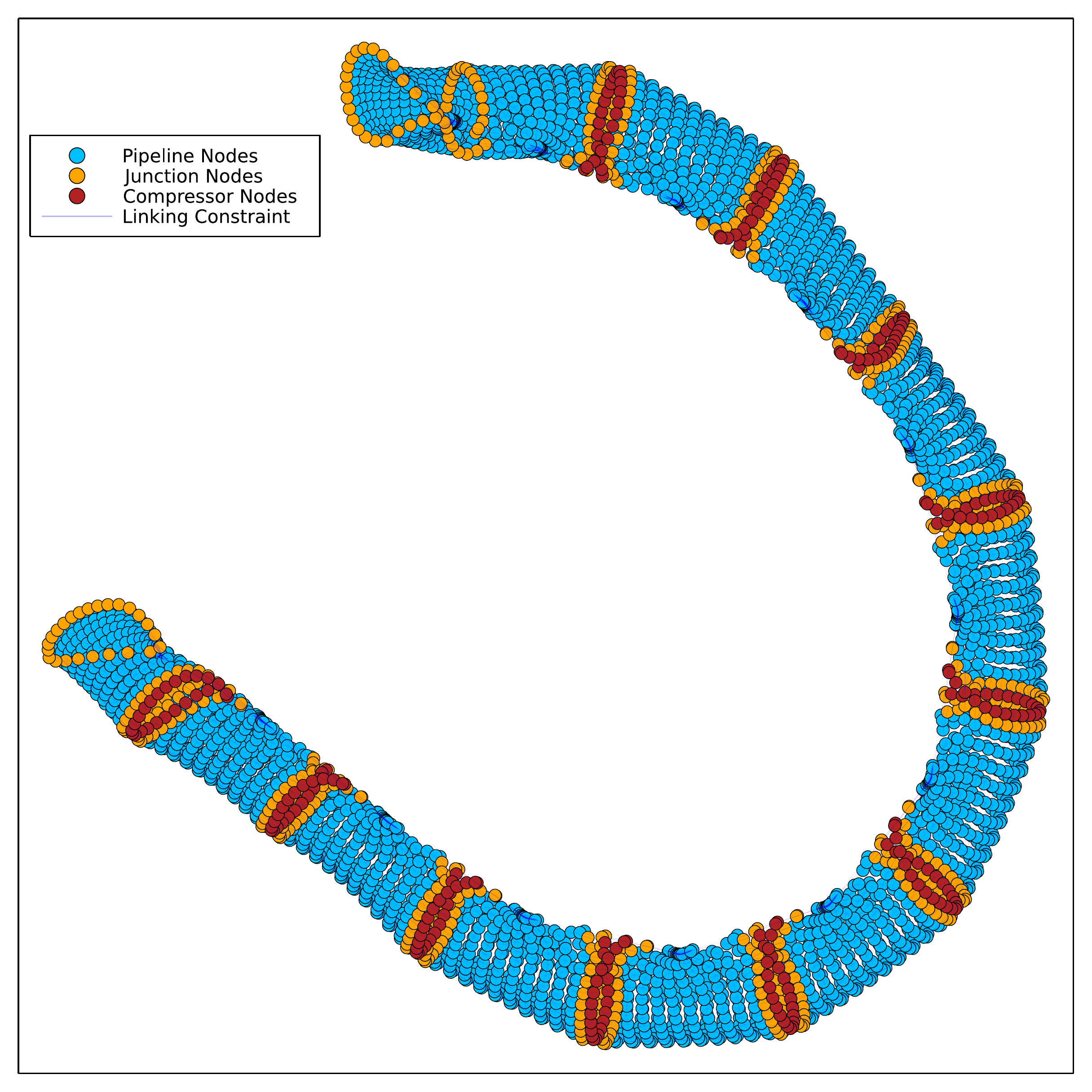}
     \caption{Visualization of a single scenario of the stochastic gas network OptiGraph}
     \label{fig:GN_one_scenario}
\end{figure}

A visualization of a scenario of this gas network, modeled as a graph, is shown in Figure \ref{fig:GN_one_scenario}. The circular nature of the structure along the time dimension arises due to periodicity constraints that require the final linepack in each pipeline to be greater than or equal to the initial linepack of that pipeline. There are 24 nodes within each ring, with one node for each time step. The spatial discretization of the pipelines is evident by the sets of ten pipeline nodes between each junction. Each node in Figure \ref{fig:GN_one_scenario} contains sets of variables and constraints, and are linked to other nodes through linking constraints.
\\ 

\begin{figure}
     \centering
     \includegraphics[scale=.5]{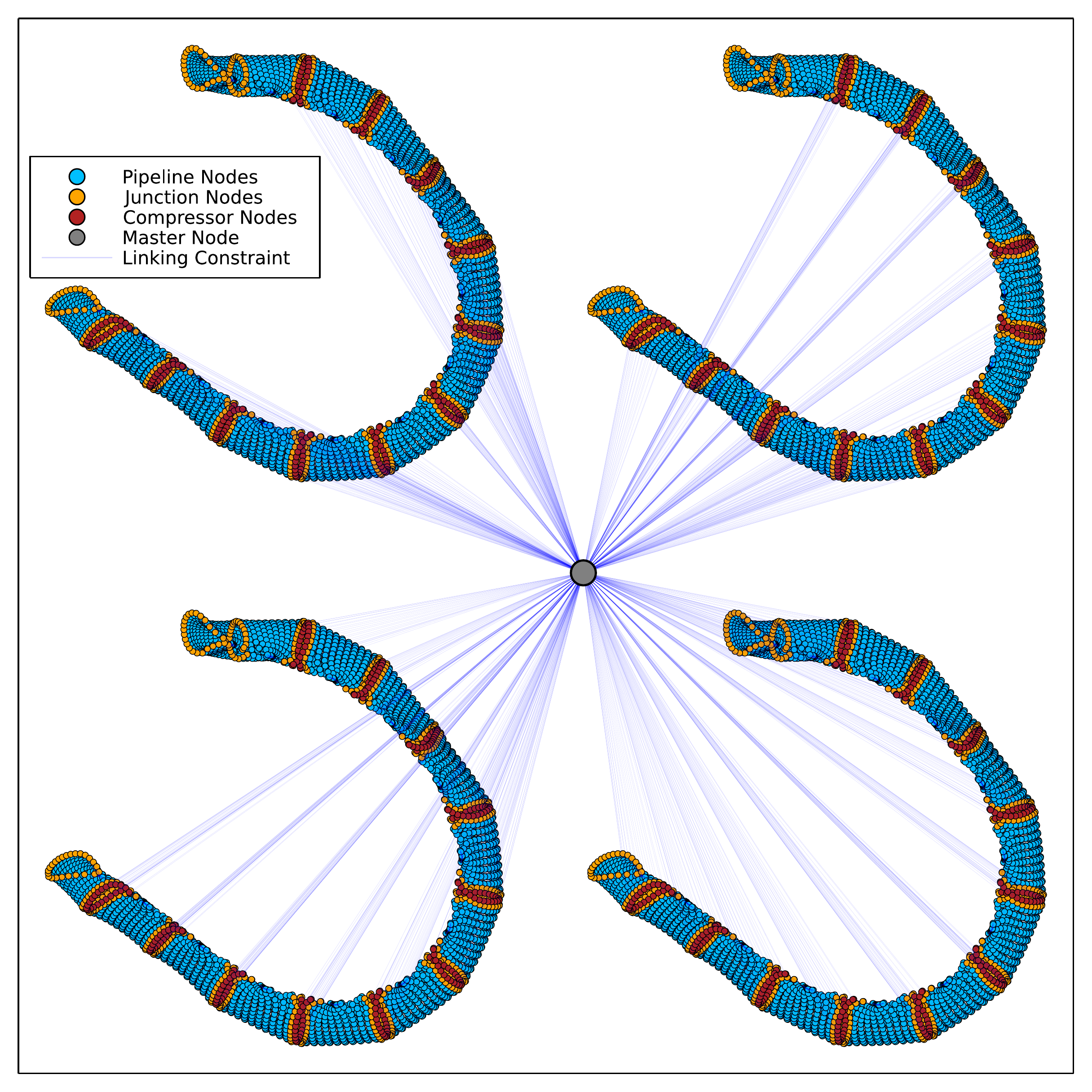}
     \caption{Visualization of a stochastic gas network OptiGraph with four scenarios. The first-stage variables on the compressor nodes are linked to a master node}
     \label{fig:GN_four_scenario}
\end{figure}

The stochastic form of this problem is visualized in Figure \ref{fig:GN_four_scenario} (with four scenarios shown). A master node in the center of the figure contains the first stage variables ($\bar{P}_{\ell, t},$ for $\ell \in \mathcal{L}_c$ and  $t \in \mathcal{T}$), and it is connected by linking constraints to each compressor node in every scenario, requiring the same compressor in each scenario to operate with the same power at each time point (see Equation \eqref{eqn:first_stage_variables}). The variables placed on this master node determine the size of the Schur complement. Because the compressors within each scenario are linked to the first-stage variables on the master node (rather than being linked directly across scenarios) the Schur complement size remains constant regardless of the number of scenarios included within the model. In the case of our formulation, {\tt MadNLP.jl} detects the two-stage structure of the problem and does not increase the Schur complement size when additional scenarios are added.
\\

Figures \ref{fig:GN_one_scenario} and \ref{fig:GN_four_scenario} show that graph-structured optimization could be useful for visualization of complex models. These visualizations also illustrate that the user can exploit at least three structures present in the problem: time, space, and scenarios. The time  structure is shown by the circular nature of the scenario graphs; the spatial structure is shown by the length of the graph from one set of junctions to the final set of junctions; and the scenario structure is shown by the separate scenario graphs that are then linked to a master node. We partitioned this problem by scenario by placing each scenario on its own subgraph and aggregating each subgraph into a node using {\tt Plasmo.jl} prior to solving with {\tt MadNLP.jl}. This resulted in the two-level tree structure seen in two-stage stochastic optimization. However, it would also be possible to construct and partition the problem by its other structures (time or space). In addition, the visualizations of other graph-based optimization models could lead to insights into additional structures and other non-obvious structures could be illuminated by using a graph-partitioning algorithm. 
\\
 
\subsection{Model Results}

\begin{figure}
  \centering
  \includegraphics[scale=.3]{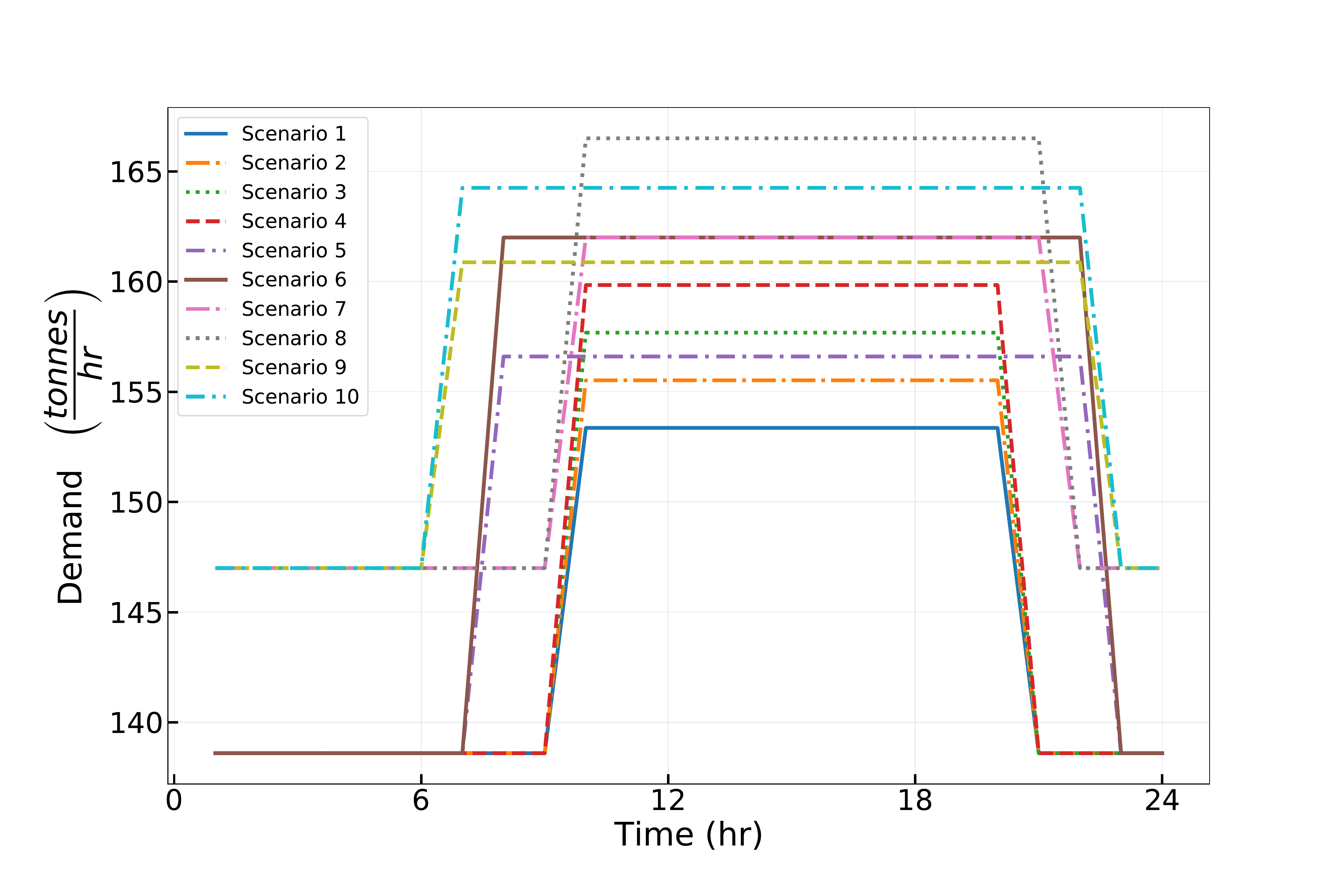}
  \caption{Natural gas demand profiles in 10 different scenarios for a stochastic natural gas network problem. The network consisted of 11 compressors, 13 pipelines, and 25 junctions. When operating with the same compressor power policy, each scenario was able to exactly meet the given demand with the exception of time point 1 of scenario 6, which over delivered by less than 2 tonnes per hour.}
  \label{fig:Gas_Demands}
\end{figure}

To demonstrate the capability of this two-stage stochastic model, we show the optimal solution of this problem under 10 scenarios for a 170 km network (pipe diameter of 0.92 m) under the formulation discussed above. The gas demand at the end of the pipeline varied across each scenario and can be seen in Figure \ref{fig:Gas_Demands}. Each scenario involved a step up and a step down, but the magnitude and location of that step could vary between scenarios. While each scenario operated under the same compressor power policy, the gas demands in each scenario were still able to be met exactly, with the exception of the first time point of scenario 6, which over delivered by less than 2 tonnes per hour. 
\\

\begin{figure}
  \centering
  \includegraphics[scale=.3]{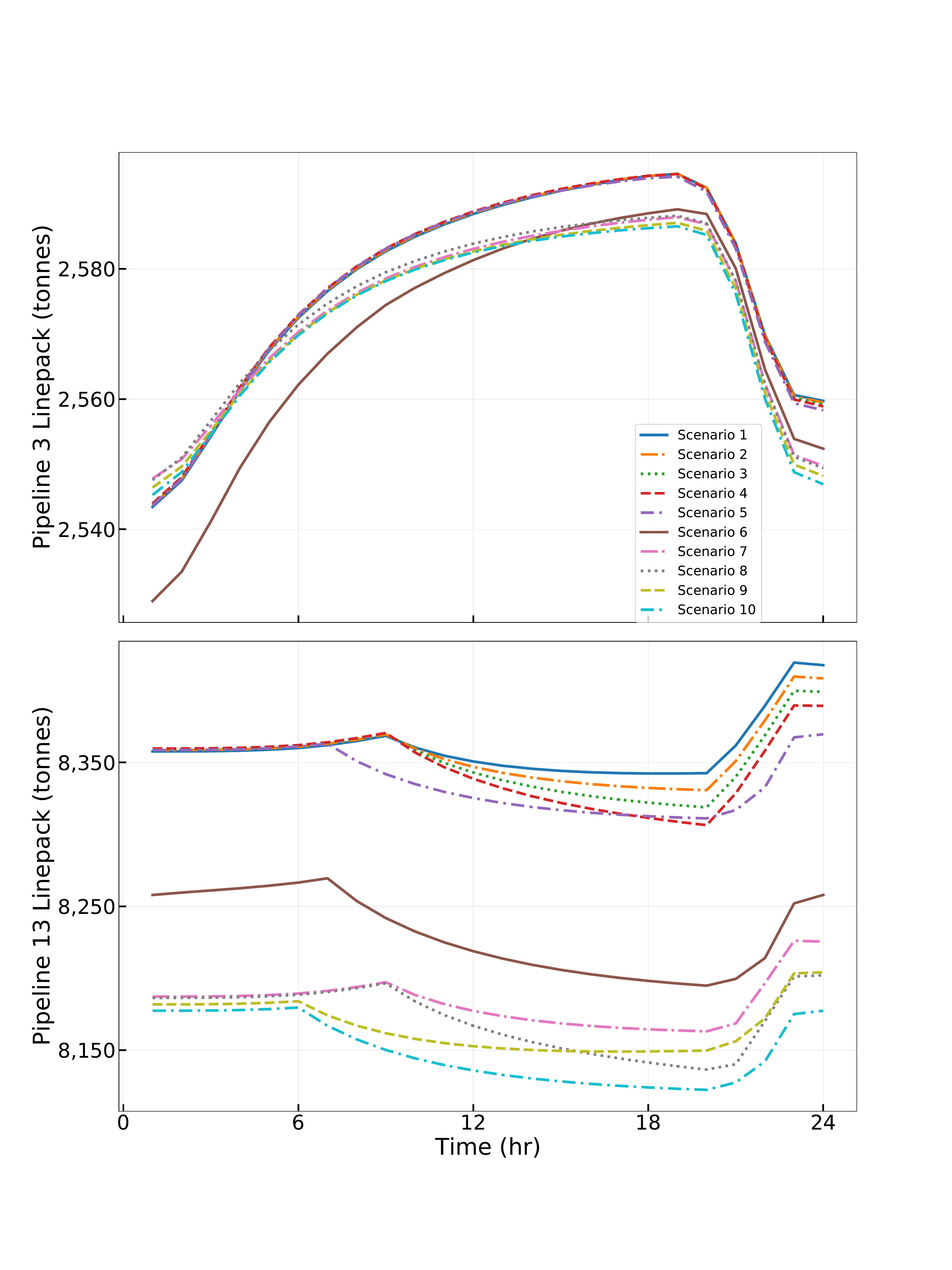}
  \caption{Linepack change over time for a stochastic natural gas network problem containing 11 compressors, 13 pipelines, and 25 junctions. The linepack in 10 different scenarios varied with time to meet varying gas demands at the end of the pipeline. Pipeline 3 (top) and pipeline 13 (bottom) are shown.}
  \label{fig:linepack_results}
\end{figure}

\begin{figure}
  \centering
  \includegraphics[scale=.3]{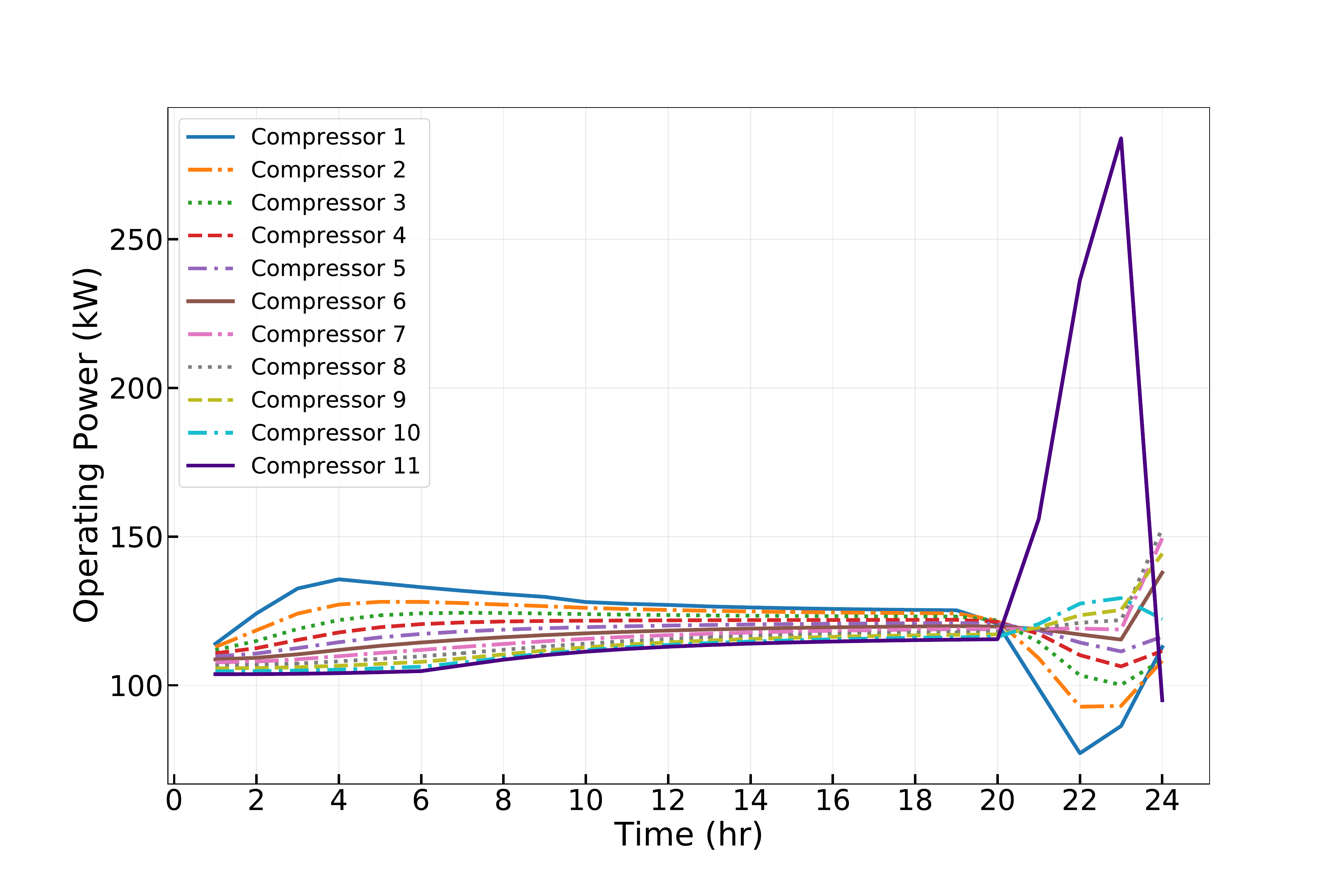}
  \caption{Compressor power policy for a stochastic natural gas network problem. The network consisted of 11 compressor, 13 pipelines, and 25 juncitons with varying gas demands at the end of the network.}
  \label{fig:Compressor_Policy}
\end{figure}

The ability to meet the varying gas demands while operating under the same compressor power policy is largely made possible by varying the pressure and linepack within each scenario. Figure \ref{fig:linepack_results} shows the linepack across different scenarios for two pipelines within the model (pipelines 3 and 13). The linepack shape over time was similar across scenarios, but the magnitude could differ to meet the required demand. The optimal compressor power policy is shown in Figure \ref{fig:Compressor_Policy}. Most compressors operated at a similar power from hour 2 until hour 20, where powers started to diverge. Compressor 11 saw a sharp increase during the last few time points. These changes in the later time periods were likely due to compressors restoring the linepack to the original values. The linepack in pipeline 13 was depleted over time until about hour 20, after which the linepack was slowly restored (Figure \ref{fig:linepack_results}). Pipeline 13 is one of two pipelines that follows compressor 11, the compressor that saw a sharp peak in power usage at hour 20.   

\subsection{Solver Performance}

\begin{figure}
  \centering
  \includegraphics[scale=.3]{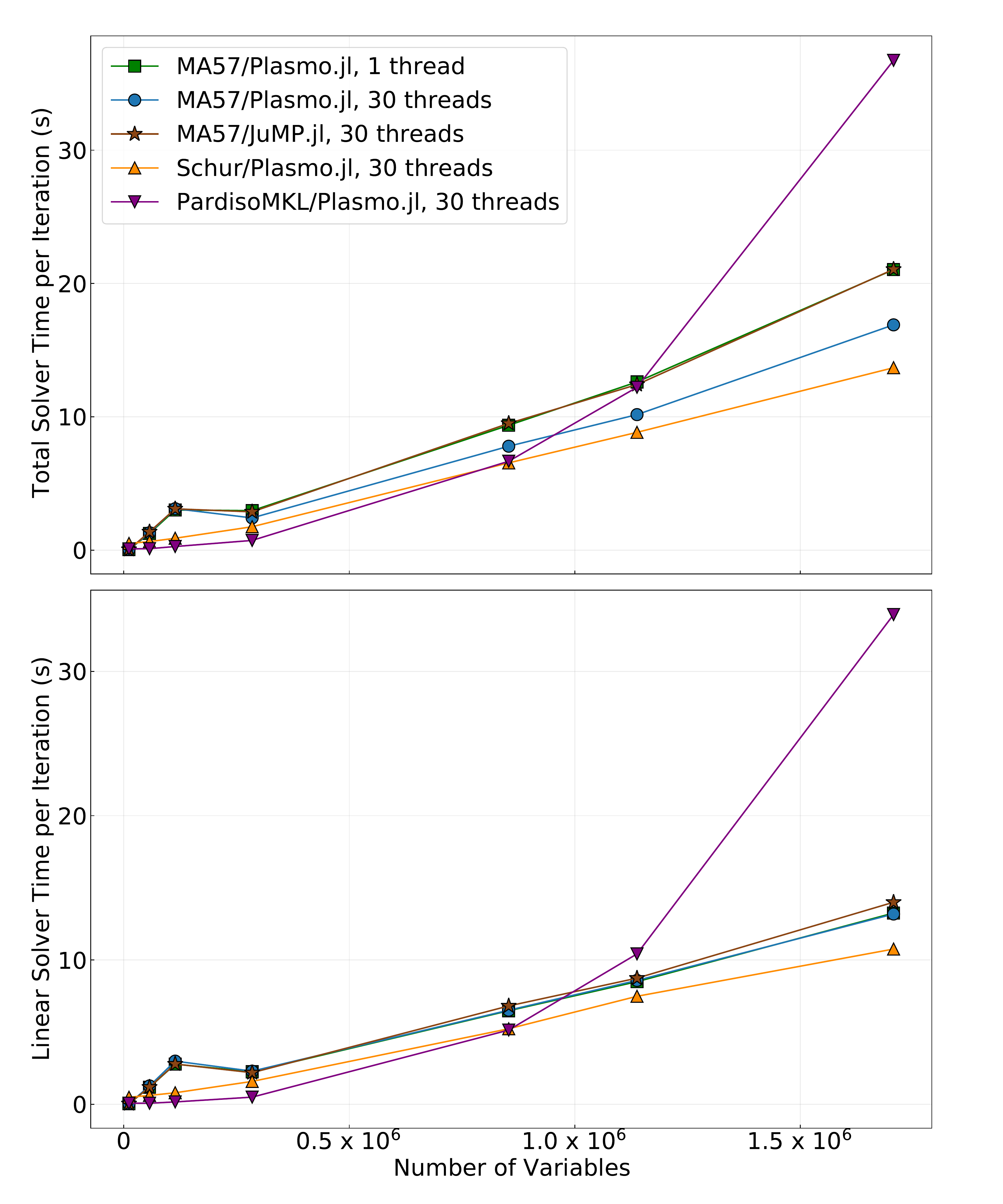}
  \caption{Total solution time per iteration (top) and linear solution time per iteration (bottom) with different problem sizes for {\tt MadNLP.jl} equipped with solvers MA57, Schur, PardisoMKL. The problem was formulated as either a {\tt Plasmo.jl} OptiGraph or a {\tt JuMP.jl} model as indicated. Values represent an average of three trials}
  \label{fig:sol_times}
\end{figure}

To test solver performance, we solved this optimization problem with different numbers of scenarios in {\tt MadNLP.jl}. Here, we used three different solver options in {\tt MadNLP.jl}: MA57, PardisoMKL, and {\tt MadNLP.jl}'s Schur solver equipped with MA57 as the linear solver. When operating with MA57, we tested {\tt MadNLP.jl} in both serial and with 30 parallel threads, and we tested PardisoMKL and Schur using 30 threads each. In each of these tests, the problem was modeled as a graph using {\tt Plasmo.jl}'s OptiGraph. For comparison to a problem not modeled as a graph, we also created a {\tt JuMP.jl} model of the system and solved that model with MA57 where {\tt MadNLP.jl} had access to 30 parallel threads (this helps compare efficiency gains from structured modeling). Each test was run on a shared-memory parallel computing server that contained a 40 core Intel(R) Xeon(R) CPU E5-2698 v4 processor running @ 2.2 GHz. Solution times for different problem sizes are shown in Figure \ref{fig:sol_times} and given in Tables \ref{tab:solution_times} and \ref{tab:solution_times_pit}. The times shown in these figures were collected by running each test four times and averaging the last three solution times. The first time was omitted to exclude compilation time. 
\\

\begin{table}[h]
  \caption{Solution times with different problem sizes for {\tt MadNLP.jl} equipped with solvers MA57, Schur, and PardisoMKL. The problem was formulated as either a {\tt Plasmo.jl} OptiGraph or a {\tt JuMP.jl} model as indicated. Values represent an average of three trials}
  \centering 
  \begin{tabular}{c rrrrrrrr} 
  \hline\hline 
   & & \multicolumn{7}{c}{Problem Size (Thousands of Variables)} \\
   Solver &  & 12 & 57 & 114 & 285 & 853 & 1,138 & 1,707 \\
  \hline
  \hline 
  \multirow{3}{*}{\shortstack{MA57/Plasmo.jl\\ 1 thread}}& iterations & 45 & 93 & 97& 129& 150& 145& 183\\ 
   & total solver time (s) & 2.23 & 116.9 & 292.6& 382.5& 1405.2& 1829.4& 3850.5\\
   & linear solver time (s) & 1.86 & 109.9 & 270.0& 291.3& 973.2& 1232.8& 2425.8\\ 
  \hline
  \multirow{3}{*}{\shortstack{MA57/Plasmo.jl\\ 30 threads}}& iterations & 45 & 93 & 97& 130 & 150& 145& 183\\ 
  & total solver time (s) & 2.42 & 121.7 & 301.7& 315.4& 1169.4& 1474.3& 3092.5\\
  & linear solver time (s) & 1.92 & 118.7 & 290.6& 299.8& 978.7& 1245.1& 2413.4\\ 
 \hline
 \multirow{3}{*}{\shortstack{MA57/JuMP.jl\\ 30 threads}}& iterations & 45 & 95 & 100 & 123 & 146 & 145 & 185 \\ 
 & total solver time (s) & 2.33 & 130.6 & 310.0& 353.7& 1388.2&1801.0 &3896.9 \\
 & linear solver time (s) & 1.86 & 113.3 &279.5 &269.9 & 995.6& 1268.0& 2590.7\\ 
\hline 
\multirow{3}{*}{\shortstack{Schur/Plasmo.jl\\ 30 threads}}& iterations & 44 & 94 & 100& 120& 139& 134& 178\\ 
& total solver time (s) & 22.3 & 60.7 & 89.6 & 211.0 & 910.1 & 1184.3 & 2436.1 \\
& linear solver time (s) & 21.6 & 58.2 & 80.2 & 190.8 & 728.0 &1003.3 & 1915.6 \\  
\hline 
\multirow{3}{*}{\shortstack{PardisoMKL/Plasmo.jl\\ 30 threads}}& iterations & 47 & 94 &98 &123 &145 &144 &177 \\
& total solver time (s) & 4.48 & 10.7 &27.0 &90.0 &967.2 &1763.7 &6512.9 \\
& linear solver time (s) & 3.68 & 7.30 &16.5 &61.4 &746.4 &1504.8 &6016.5 \\
\hline 
  \end{tabular}
  \label{tab:solution_times}
\end{table}

\begin{table}[h]
  \caption{Solution times per iteration with different problem sizes for {\tt MadNLP.jl} equipped with solvers MA57, Schur, and PardisoMKL. The problem was formulated as either a {\tt Plasmo.jl} OptiGraph or a {\tt JuMP.jl} model as indicated} 
  \centering
  \begin{tabular}{c rrrrrrrr}
  \hline
  \hline 
   & & \multicolumn{7}{c}{Problem Size (Thousands of Variables)} \\
   Solver & Value per Iteration  & 12 & 57 & 114 & 285 & 853 & 1,138 & 1,707 \\
  \hline
  \hline  
  \multirow{2}{*}{\shortstack{MA57/Plasmo.jl\\ 1 thread}}& total solver time (s) & 0.05 & 1.26 & 3.02 & 2.97 & 9.37&12.62 &21.04 \\ % 
   & linear solver time (s) & 0.04& 1.18 &2.78 &2.26 &6.49 &8.50 &13.26 \\ 
  \hline 
  \multirow{2}{*}{\shortstack{MA57/Plasmo.jl\\ 30 threads}}& total solver time (s) & 0.05 &1.31 &3.11 &2.42 &7.80 &10.17 &16.90\\ % 
   & linear solver time (s) & 0.04&1.28  &3.00 &2.30 &6.52 &8.59 &13.19 \\
  \hline 
  \multirow{2}{*}{\shortstack{MA57/JuMP.jl\\ 30 threads}}& total solver time (s) &0.05 &1.37 &3.10 &2.88 &9.51 &12.42 &21.06\\ % 
   & linear solver time (s) &0.04 &1.19  &2.79 &2.19 &6.82 &8.75 &14.00 \\ 
  \hline
  \multirow{2}{*}{\shortstack{Schur/Plasmo.jl\\ 30 threads}}& total solver time (s) & 0.51& 0.65& 0.90&1.76 &6.54 &8.84 &13.69\\ % 
   & linear solver time (s) & 0.49& 0.62 &0.80 &1.59 &5.24 &7.49 &10.76 \\ 
  \hline
  \multirow{2}{*}{\shortstack{PardisoMKL/Plasmo.jl\\ 30 threads}}& total solver time (s) &0.10 &0.11 &0.27 &0.73 &6.67 &12.22 &36.74\\ % 
  & linear solver time (s) & 0.08& 0.08 &0.17 &0.50 &5.15 &10.42 &33.94 \\ 
  \hline 
\hline
  \end{tabular}
  \label{tab:solution_times_pit}
\end{table}

The results in Figure \ref{fig:sol_times} and Tables \ref{tab:solution_times} and \ref{tab:solution_times_pit} provide comparisons between different solver options. At small problem sizes ($\sim$ 10,000 variables), {\tt MadNLP.jl} with MA57 outperformed {\tt MadNLP.jl} with Schur or PardisoMKL based on total time per iteration. However, solver times saw the more drastic differences as the number of variables increased. MA57 with access to one thread took the longest time to run per iteration from about 100,000 to 1.1 million variables, but it was surpassed by PardisoMKL when run with 1.7 million variables. {\tt MadNLP.jl} operating with MA57 with access to 30 threads was 20\% faster per iteration at 1.7 million variables than when {\tt MadNLP.jl} had access to only one thread. The Schur solver did better than any run of MA57 after the problem was scaled beyond $\sim$ 10,000-50,000 variables. For the largest problem sizes, the Schur solver was 19\% faster per iteration than {\tt MadNLP.jl} operating with MA57 and 30 threads, and it was 35\% faster per iteration than {\tt MadNLP.jl} operating with MA57 and one thread . PardisoMKL also performed well at lower problem sizes ($\le$ $\sim$600,000 variables), but showed a substantial increase in solver time as problem sizes got very large, ultimately taking 2.7 times the time per iteration compared to the Schur solver at 1.7 million variables. 
\\

Notably, there were significant differences between the total time per iteration at the largest problem size for {\tt MadNLP.jl} operating in parallel or serial on either {\tt Plasmo.jl}'s OptiGraph or {\tt JuMP.jl}'s model. Solving the {\tt JuMP.jl} model with 30 threads took the same amount of time as solving a {\tt Plasmo.jl} OptiGraph with 1 thread while operating with MA57. This is because {\tt MadNLP.jl} does not parallelize the {\tt JuMP.jl} model since it is not constructed as a graph. However, when solving a {\tt Plasmo.jl} OptiGraph with 30 threads rather than 1 thread, {\tt MadNLP.jl} operating with MA57 had a reduction of more than four seconds in total solver time per iteration even though the linear solver times per iteration were virtually identical. This is because, while MA57 has limited parallelization capability (evidenced by the linear solver results), {\tt MadNLP.jl} is able to parallelize function and derivative evaluations outside of the linear solver to reduce total solution time. This highlights that non-trivial improvements in computational time can be achieved by exploiting structure outside the linear solver. 
\\

These results highlight some of the benefits of graph-structured optimization and the capabilities of {\tt MadNLP.jl}. First, graph-structured optimization and {\tt MadNLP.jl} allow for function and derivative evaluations to be run in parallel, even when the linear solver itself has limited parallelization abilities. In addition, graph-structured optimization provides an efficient method for implementing decomposition schemes such as Schur decomposition that can be run in parallel. While the Schur solver was equipped with MA57 as a linear solver, the linear solver time per iteration was less than MA57 alone beyond problem sizes of $\sim$ 10,000-50,000 variables. Thus graph-structured optimization has the ability to reduce solution times for complex problems. {\tt MadNLP.jl} takes advantage of this ability and efficiently implements these concepts. 

%%%%%%%%%%%%%%%%%%%%%%%%%%%%%%%%%%%%%%%%%%

\section{Conclusions and Future Work}

We highlighted many of the benefits of graph-structured optimization, including visualization capabilities, partitioning schemes, the parallelization of function and derivative evaluations, and decomposition schemes to reduce solver time. These advantages are utilized by {\tt Plasmo.jl} and {\tt MadNLP.jl} to model and solve graph-structured optimization problems. We gave an illustrative example of the visualization and partitioning capabilities of {\tt Plasmo.jl}, and we applied {\tt MadNLP.jl} to solve a case study of a stochastic natural gas pipeline using multiple solvers available within {\tt MadNLP.jl}. The results of this case study suggest that significant time can be saved by parallelizing function and derivative evaluations and using decomposition schemes to parallelize the linear solver. 
\\

While graph-based modeling offers many benefits, this paradigm can be expanded in a number of ways. We are interested in increasing the ability of {\tt MadNLP.jl} to operate with subgraphs without aggregating some subgraphs into nodes prior to solution. In addition, we are interested in additional decomposition schemes that could also be applied to graph-structured optimization.  Additional work can be done to test the capabilities and applications of Schwarz decomposition, particularly its ability to operate as an iterative solver.  

%%%%%%%%%%%%%%%%%%%%%%%%%%%%%%%%%%%%%%%%%%

\section*{Acknowledgments}
This work was partially supported by the U.S. Department of Energy under grant DE-SC0014114. We also acknowledge the contributions of Jordan Jalving for his recommendations for modeling and interfacing with {\tt Plasmo.jl}. 

%%%%%%%%%%%%%%%%%%%%%%%%%%%%%%%%%%%%%%%%%%
\bibliography{MadNLP_Schur_Paper}

\end{document}